\documentclass[journal]{IEEEtran}

\usepackage{mathtools}
\usepackage{amscd}
\usepackage{amsthm,amssymb}
\usepackage{amsmath}
\usepackage{cite}
\usepackage[utf8]{inputenc} 
\usepackage{calc,ifthen,moreverb}
\usepackage{graphicx,color,psfrag}
\usepackage{lscape,pdflscape}
\usepackage{multicol}
\usepackage{dutchcal}
\usepackage{url}
\usepackage{soul}
\usepackage[caption=false,font=footnotesize]{subfig}
\usepackage{framed}
\usepackage{slashed}
\usepackage{bm}
\usepackage[table,xcdraw]{xcolor}
\usepackage{relsize}
\usepackage{cancel}
\usepackage{cleveref}
\usepackage{adjustbox} 
\usepackage{mdframed}

 \def\0z{{\boldsymbol{0}}}
 \def\1z{{\boldsymbol{1}}}
 \def\2z{{\boldsymbol{2}}}
 \def\3z{{\boldsymbol{3}}}
 \def\4z{{\boldsymbol{4}}}
 \def\5z{{\boldsymbol{5}}}
 \def\6z{{\boldsymbol{6}}}
 \def\7z{\boldsymbol{7}}
 \def\8z{\boldsymbol{8}}
 \def\9z{\boldsymbol{9}}

\usepackage{pifont}

\newcommand{\eq}{\begin{equation}}
\newcommand{\qe}{\end{equation}}
\newcommand{\eqarray}{\begin{array}{ll}}
\newcommand{\qearray}{\end{array}}

\newcommand{\norm}[1]{\left\lVert #1 \right\rVert} 

\def\*#1{\bm{#1}} 
\def\t~#1{\widetilde{#1}} 
\def\b~#1{\bm{\widetilde{#1}}} 
\def\expec{\mathop{{}\mathbb{E}}}
\def\rt#1{{#1}^{*}} 
\def\rtb#1{\bm{#1}^{*}} 
\def\iprod{ } 
\def\basis{\lambda} 

\theoremstyle{definition} 
\theoremstyle{definition} 

\begin{document}
%
\title{Geometric-Algebra Adaptive Filters}
%
%
%

\author{\IEEEauthorblockN{Wilder~B.~Lopes\IEEEauthorrefmark{1},~\IEEEmembership{Member,~IEEE,}
        Cassio~G.~Lopes\IEEEauthorrefmark{2},~\IEEEmembership{Senior~Member,~IEEE}}
        \thanks{\IEEEauthorrefmark{1}R\&D Dept. at UCit (ucit.fr) and OpenGA.org, wil@openga.org. \IEEEauthorrefmark{2}Dept. of Electronic Systems Engineering, University of Sao Paulo, Brazil, cassio@lps.usp.br.}}%
    

%
%

\markboth{}%
{Shell \MakeLowercase{\textit{et al.}}: Bare Demo of IEEEtran.cls for Journals}
%

\maketitle

\begin{abstract}
	This paper presents a new class of adaptive filters, namely Geometric-Algebra Adaptive Filters (GAAFs). They are generated by formulating the underlying minimization problem (a deterministic cost function) from the perspective of Geometric Algebra (GA), a comprehensive mathematical language well-suited for the description of geometric transformations. Also, differently from standard adaptive-filtering theory, Geometric Calculus (the extension of GA to differential calculus) allows for applying the same derivation techniques regardless of the type (subalgebra) of the data, i.e., real, complex numbers, quaternions, etc. Relying on those characteristics (among others), a deterministic quadratic cost function is posed, from which the GAAFs are devised, providing a generalization of regular adaptive filters to subalgebras of GA. From the obtained update rule, it is shown how to recover the following least-mean squares (LMS) adaptive filter variants: real-entries LMS, complex LMS, and quaternions LMS. Mean-square analysis and simulations in a system identification scenario are provided, showing very good agreement for different levels of measurement noise.
\end{abstract}


\ifCLASSOPTIONpeerreview
\begin{center} \bfseries EDICS Category: ASP-ANAL \end{center}
\fi
%
\IEEEpeerreviewmaketitle

\section{Introduction}\label{cap:intro}
\begin{figure}[!t]
	\centering
	\includegraphics[angle=-90, width=0.30\textwidth]{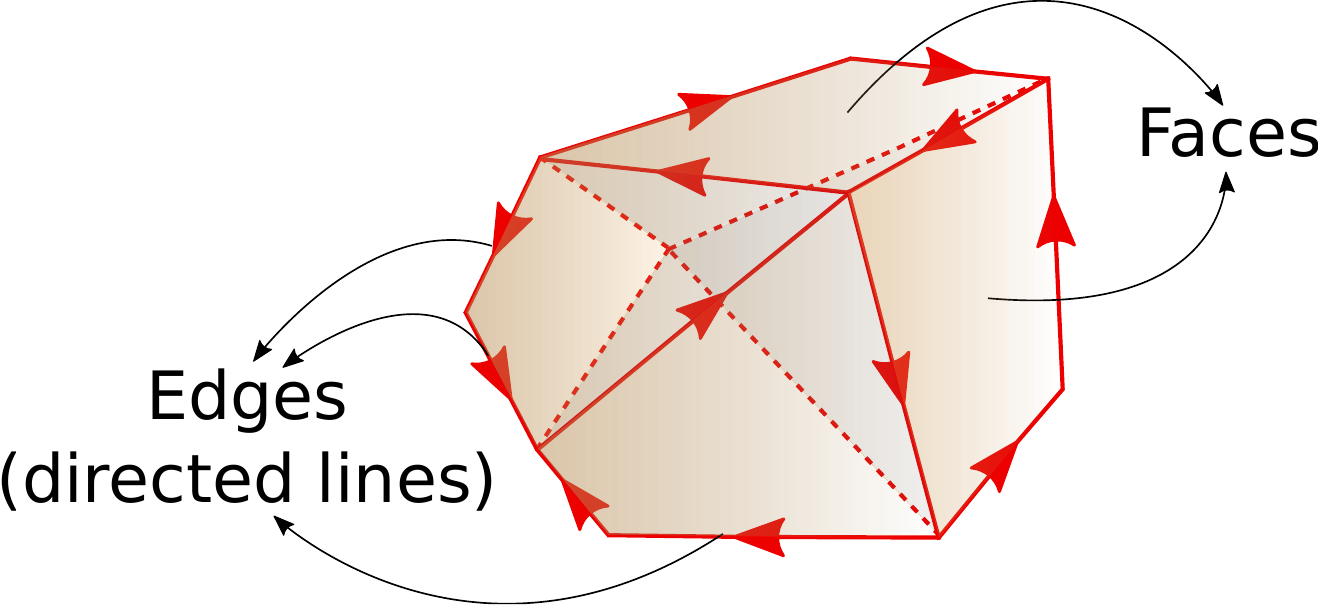}
	\caption{A polyhedron (3-dimensional polytope) can be completely described by the geometric multiplication of its edges (oriented lines, vectors), which generate the faces and hypersurfaces (in the case of a general n-dimensional polytope).}
	\label{fig:polytope}
\end{figure}

\IEEEPARstart{A}{daptive} filters (AFs), usually described via linear algebra (LA) and standard vector calculus~\cite{Sayed08}, can be interpreted as instances for geometric estimation, since the weight vector to be estimated represents a directed line in an underlying vector space, i.e, a n-dimensional vector encodes the length and direction of an edge in a n-dimensional polytope (see \figurename~\ref{fig:polytope}). However, to estimate the polytope as a whole (edges, faces, and so on), a regular AF designed in light of LA might not be very helpful. 

Indeed, LA has limitations regarding the representation of geometric structures~\cite[p. 20]{hestenes1999newfoundations}. Take for instance the matrix product between two vectors: it always results in either a scalar or a matrix (2-dimensional array of numbers). Thus, one may wonder if it is possible to construct a new kind of product that takes two vectors (directed lines or edges) and returns a hypersurface (the face of an n-dimensional polytope); or takes a vector and a hypersurface and returns another polytope. Similar ideas have been present since the advent of algebra, in an attempt to establish a deep connection with geometry~\cite{crowe1967history,kleiner2007history}.

The \emph{geometric product}, which is the product operation of GA~\cite{hestenes1987GAcalculus,Hitzer_Introduction}, captures the aforementioned idea. It allows one to map a set of vectors not only onto scalars, but also onto hypersurfaces and n-dimensional polytopes. The use of GA increases the portfolio of geometric shapes and transformations one can represent. Also, its extension to calculus, namely \emph{geometric calculus} (GC), allows for a clear and compact way to perform calculus with hypercomplex quantities, i.e., elements that generalize complex numbers for higher dimensions~\cite{crowe1967history,kleiner2007history,hestenes1987GAcalculus,hestenes1999newfoundations,DoranDissertation,Hitzer_Introduction,doran2003geometric,2016introductionVaz,Dorst2007GACV}.     

GA-based AFs were first introduced in~\cite{SPL,WACV}, where they were successfully employed to estimate the geometric transformation (rotation and translation) that aligns a pair of 3D point clouds (a recurrent problem in computer vision), while keeping a low computational footprint. This work takes a step further: it uses GA and GC to generate a new class of AFs able to naturally estimate hypersurfaces, hypervolumes, and elements of greater dimensions (\emph{multivectors}), generalizing regular AFs from the literature.  Filters like the regular least-mean squares (LMS) -- real entries, the Complex LMS (CLMS) -- complex entries~\cite{Widrow1975}, and the Quaternion LMS (QLMS) -- quaternion entries~\cite{Took_QLMS,Mandic2011,Jahanchahi2012}, are recovered as special cases of the more comprehensive GA-LMS introduced herein.  

The text is organized as follows. Section~\ref{cap:fundamentals_of_GA} covers the transition from linear to geometric algebra, presenting the fundamentals of GA and GC. In Section~\ref{cap:linear_estimation_GA}, the standard quadratic cost function, usually adopted in adaptive filtering, is recovered as a particular case of a more comprehensive quadratic cost function that can only be written using the geometric product. In Section~\ref{cap:GAAFs_standard} the gradient of that cost function is calculated via GC techniques and the GA-LMS is formulated. Section~\ref{sec:mean-square-analysis} provides a mean-square analysis (steady state) with the support of the energy conservation relations~\cite{Sayed08}. Experiments are shown in Section~\ref{cap:applications_GAAFs} to assess the performance of the GAAFs against the theoretical analysis in a system-identification scenario. Finally, conclusion remarks are presented in Section~\ref{cap:conclusion}. 

\subsubsection*{Remarks about notation} Table~\ref{table:notation} summarizes the notation convention. While those are deterministic quantities, their boldface versions represent random quantities. The name \emph{array of multivectors} was chosen to avoid confusion with vectors in $\mathbb{R}^n$, which in this text have the usual meaning of \emph{collection of real numbers} (row or column). In this sense, an array of multivectors can be interpreted as a ``vector'' that allows for hypercomplex entries (numbers constructed by adding ``imaginary units'' to real numbers~\cite{kantor2011hypercomplex}).

\begin{table}[t!]
	\centering
	\caption{Summary of notation.}
	\label{table:notation}
	\begin{tabular}{l|l}
		Symbol                 & Definition                                                                                                                               \\ \hline
		a                      & Arrays of multivectors, vectors and scalars.                                                                                                                      \\
		A                      & General multivector or matrix.                                                                                    \\
		r                      & Rotor (type of multivector).                                                                                                              \\
		$a(i), a_i, A(i)$ & \begin{tabular}[c]{@{}l@{}}Time-varying scalar, vector (or array of multivectors),\\ and general multivector (or matrix), respectively.\end{tabular}
	\end{tabular}
\end{table}


\section{From Linear to Geometric Algebra}
\label{cap:fundamentals_of_GA}

To start transitioning from LA to GA, one needs to recall the definition of an algebra over the reals~\cite{vinberg2003course,Hitzer_Introduction,2016introductionVaz}: a vector space $\mathcal{V}$ over the field $\mathbb{R}$, equipped with a bilinear map $\mathcal{V}\times\mathcal{V} \rightarrow \mathcal{V}$ denoted by $\circ$ (the \textit{product operation} of the algebra), is said to be an \emph{algebra} over $\mathbb{R}$ if the following relations hold $\forall \{a, b, c\} \in \mathcal{V}$ and $\{\alpha, \beta\} \in \mathbb{R}$,
\begin{equation} 
\eqarray 
(a + b)\circ c &= a\circ c + b\circ c \text{ (Left distributivity)}, \\
c\circ(a + b) &= c\circ a + c\circ b \text{ (Right distributivity)}, \\
(\alpha a)\circ(\beta b) &= (\alpha\beta)(a\circ b) \text{ (Compatibility with scalars).}
\qearray
\label{def:definition_algebra}
\end{equation}

Linear (matrix) algebra, utilized to describe adaptive filtering theory, is constructed from the definition above. The elements of this algebra are matrices and vectors, which multiplied among themselves via the \textit{matrix product} generate new matrices and vectors. Additionally, to express the notion of vector length and angle between vectors, LA adopts the bilinear form $\mathcal{V}\times\mathcal{V} \rightarrow \mathbb{R}$, i.e., \textit{inner product}, which returns a real number as a result of the multiplication between two vectors in $\mathcal{V}$ (one says that $\mathcal{V}$ is a normed vector space)~\cite[p. 180]{vinberg2003course}.    

Geometric (Clifford) Algebra derives from~\eqref{def:definition_algebra}, however with a different product operation. Such product, called \textit{geometric product}, is what allows for GA to be a mathematical language that unifies different algebraic systems trying to express geometric relations/transformations, e.g., rotation and translation. The following systems are examples of algebras integrated into GA theory: vector/matrix algebra, complex numbers, and quaternions~\cite{hestenes1987GAcalculus,hestenes1999newfoundations,Hitzer_Introduction,2016introductionVaz}. Such unifying quality is put into use in this work to expand the capabilities of AFs.  

The fundamentals of GA are provided in the sequel. For an in-depth discussion, the reader is referred to~\cite{hestenes1987GAcalculus,hestenes1999newfoundations,DoranDissertation,Hitzer_Introduction,doran2003geometric,2016introductionVaz,Dorst2007GACV,Lasenby1998,perwass2009geometric,hildenbrand2012foundations,dorst2012applications,sommer2013geometric,crowe1967history,kleiner2007history}.  

\subsection{Fundamentals of Geometric Algebra}
\label{sec:GA_from_VecSpaces}

Consider $\{a,b\}$ vectors in $\mathbb{R}^n$, i.e., \emph{arrays with real entries}. The inner product $a \cdot b$ is the standard bilinear form that describes vector length and angle between vectors in linear (matrix) algebra. This way, $a\cdot b$ results in a scalar,
\vspace*{-2mm} 
\eq 
a \cdot b = |a||b|cos\theta,
\label{def:inner_product_vectors}
\vspace*{-1.5mm}
\qe 
in which $\theta$ is the angle between $a$ and $b$, and $|\cdot|$ denotes the vector magnitude (norm). Additionally, the inner product is commutative, i.e., $a \cdot b = b \cdot a$.

\begin{figure}[!t]
	\centering
	\includegraphics[width=0.45\textwidth]{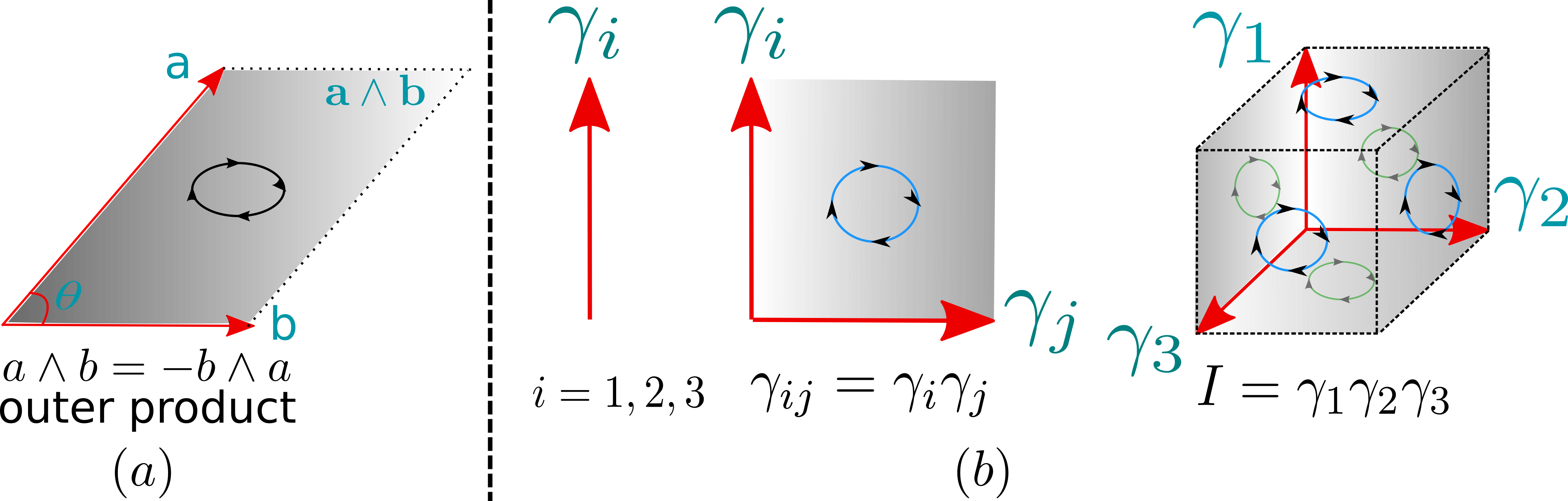}
	\caption[Outer product and basis of $\mathcal{G}(\mathbb{R}^3)$]{(a) Visualization of the outer product in $\mathbb{R}^3$. The orientation of the circle defines the orientation of the surface (bivector). (b) The elements of $\mathcal{G}(\mathbb{R}^3)$ basis (besides the scalar $1$): $3$ vectors, $3$ bivectors (oriented surfaces) $\gamma_{ij}$, and the trivector $I$ (pseudoscalar/oriented volume). Note that in $\mathbb{R}^3$, $\langle A \rangle_g = 0$, $g > 3$~\cite[p.42]{hestenes1999newfoundations}.}
	\label{fig:multivectors_in_R3}
\end{figure}

The outer product $a \wedge b$ is the usual product of the exterior algebra introduced by Grassmann's \textit{Ausdehungslehre} (theory of extension)~\cite{crowe1967history,grassmann2000ausdehnungslehre,hestenes1999newfoundations,2016introductionVaz}. It captures the geometric fact that two nonparallel directed segments determine a parallelogram, a notion which can not be described by the inner product. The multiplication $a\wedge b$ results in an \emph{oriented surface} or \emph{bivector} (see \figurename~\ref{fig:multivectors_in_R3}a). Such a surface can be interpreted as the parallelogram (hyperplane) generated when vector $a$ is swept on the direction determined by vector $b$. Alternatively, the outer product can be defined as a function of the angle $\theta$ between $a$ and $b$ 
\vspace*{-2mm}
\eq 
a \wedge b =  I_{a,b}|a||b|sin\theta,
\label{def:outer_product_vectors}
\vspace*{-1.5mm}
\qe 
where $I_{a,b}$ is the \emph{unit bivector}\footnote{An unit bivector is the result of the outer product between two unit vectors, i.e., vectors with unitary norm.} that defines the orientation of the hyperplane $a \wedge b$~\cite[p.66]{hestenes1999newfoundations}. Note that in the particular case of 3D Euclidean space, \eqref{def:outer_product_vectors} is reduced to the \textit{cross product} $a \times b = p|a||b|sin\theta$, where $p$ is the unit vector normal to the plane containing $\{a, b\}$. From \figurename~\ref{fig:multivectors_in_R3}a it can be concluded that the outer product is noncommutative, i.e., $a \wedge b = -b \wedge a$: the orientation of the surface generated by sweeping $a$ along $b$ is opposite to the orientation of the surface generated by sweeping $b$ along $a$.

Finally, the geometric product\footnote{In this text, from now on, all products are geometric products, unless otherwise noted.} of vectors $a$ and $b$ is denoted by their juxtaposition $ab$ and defined as 
\vspace*{-2mm}
\eq 
ab \triangleq a \cdot b + a \wedge b,
\label{def:geometric_product_vectors}
\vspace*{-1.5mm}
\qe
in terms of the inner ($\cdot$) and outer ($\wedge$) products~\cite[Sec. $2.2$]{Hitzer_Introduction}.
Note that due to $a \wedge b = -(b \wedge a)$ the geometric product is \emph{noncommutative}, resulting in $ab = -ba$, and it is \emph{associative}, $a(bc) = (ab)c$, $\{a,b,c\} \in \mathbb{R}^n$.

In linear algebra the fundamental elements are matrices/vectors. In a similar way, in GA the basic elements are the so-called \textit{multivectors} (Clifford numbers). The structure of a multivector can be seen as a generalization of complex numbers and quaternions for higher dimensions. For instance, a complex number $\alpha + j\beta$ has a scalar part $\alpha$ combined with an imaginary part $j\beta$; quaternions like  $\alpha + i\beta_1 + j\beta_2 + k\beta_3$ expand that by adding two extra imaginary-valued parts. Multivectors generalize this structure, which contains one scalar part and several other parts named \textit{grades}. Thus, a general multivector $A$ has the form
\vspace*{-2mm}
\begin{equation}
A = \langle A \rangle_0 + \langle A \rangle_1 + \langle A \rangle_2 + \cdots = \sum_{g} \langle A \rangle_g \text{,} 
\label{def:multivector}
\vspace*{-1.5mm}
\end{equation} 
which is comprised of its g-grades (or g-vectors) $\langle\cdot\rangle_g$, e.g., $g=0$ (scalars), $g=1$ (vectors), $g=2$ (bivectors, generated via the geometric multiplication of two vectors), $g=3$ (trivectors, generated via the geometric multiplication of three vectors), and so on. The grade $g=0$ (scalar) is also denoted as $\langle A \rangle_0 \equiv \langle A \rangle$. Additionally, in $\mathbb{R}^n$, $\langle A \rangle_g = 0$, $g > n$~\cite{hestenes1999newfoundations}. The ability to group together scalars, vectors, and hypercomplex quantities in a unique element (multivector) is the foundation on top of which GA theory is built -- it allows for ``summing apples and oranges'' in a well-defined fashion. Section~\ref{sec:subalgebras_isomorphism} will show how to recover complex numbers and quaternions as special cases of~\eqref{def:multivector}.

The multivectors that form the basis of a geometric algebra over the vector space $\mathcal{V}$, denoted $\mathcal{G}(\mathcal{V})$, are obtained by multiplying the $n$ vectors that compose the orthonormal basis of $\mathcal{V}$ via the geometric product~\eqref{def:geometric_product_vectors}. This action generates $2^n$ multivectors, which implies that $\dim\{\mathcal{G}(\mathcal{V})\} = 2^n$~\cite[p. 19]{hestenes1987GAcalculus}. Those $2^n$ multivectors are called \emph{blades} of $\mathcal{G}(\mathcal{V})$. For the special case $n = 3 \Rightarrow \mathcal{V} = \mathbb{R}^3$, with orthonormal basis $\{\gamma_1,\gamma_2,\gamma_3\}$, the procedure above yields the following blades
\vspace*{-2mm}
\eq 
\{1,\gamma_1,\gamma_2,\gamma_3,\gamma_{12},\gamma_{23},\gamma_{31},I\}, 
\label{eq:orthonormal_basis_R3}
\vspace*{-1.5mm}
\qe 
which together are a basis for $\mathcal{G}(\mathbb{R}^3)$ with dimension $2^3 = 8$. Note that~\eqref{eq:orthonormal_basis_R3} has one scalar, three orthonormal vectors $\gamma_{i}$ (basis for $\mathbb{R}^3$), three bivectors (oriented surfaces) $\gamma_{ij} \triangleq \gamma_{i}\gamma_{j}=\gamma_{i}\wedge\gamma_{j}, i \neq j$ ($\gamma_{i}\cdot\gamma_{j} = 0, i \neq j$), and one trivector (\emph{pseudoscalar}) $I \triangleq \gamma_1\gamma_2\gamma_3 = \gamma_{123}$ (\figurename~\ref{fig:multivectors_in_R3}b). In general, the pseudoscalar  $I$ is defined as the multivector of highest grade in an algebra $\mathcal{G(V)}$. It commutes with any multivector in $\mathcal{G(V)}$, hence the name pseudoscalar~\cite[p. 17]{hestenes1987GAcalculus}.    

To illustrate the geometric multiplication between elements of $\mathcal{G}(\mathbb{R}^3)$, take two multivectors $C = \gamma_1$ and $D = 2\gamma_1 + 4\gamma_3$. Then, $CD = \gamma_1(2\gamma_1 + 4\gamma_3) = \gamma_1\cdot(2\gamma_1 + 4\gamma_3) + \gamma_1\wedge(2\gamma_1 + 4\gamma_3) = 2 + 4(\gamma_1\wedge\gamma_3) = 2 + 4\gamma_{13}$ (a scalar plus a bivector). Another example is provided to highlight the \emph{reverse} of a multivector $A$, which is the GA counterpart of complex conjugation in linear algebra, defined as
\vspace*{-2mm}
\begin{equation}
\widetilde{A} \triangleq \sum_{g=0}^{n} (-1)^{g(g-1)/2}\langle A \rangle_g\text{.}
\label{eq:reversion}
\vspace*{-2mm}
\end{equation}
Thus, given $A = \langle A \rangle_0 + \langle A \rangle_1 + \langle A \rangle_2$, its reverse is $\widetilde{A} = \widetilde{\langle A \rangle_0} + \widetilde{\langle A \rangle_1} + \widetilde{\langle A \rangle_2} = \langle A \rangle_0 + \langle A \rangle_1 - \langle A \rangle_2$. Therefore, the reverse of the $CD$ multiplication above is $\widetilde{(CD)} = 2 - 4\gamma_{13}$. Note that since the $0$-grade of a multivector is not affected by reversion, mutually reverse multivectors, say $A$ and $\t~A$, have the same $0$-grade, $\langle A \rangle_0 = \langle \widetilde{A} \rangle_0$.

\subsection{Subalgebras and Isomorphisms}
\label{sec:subalgebras_isomorphism}
The GAAFs are designed to compute with any multivector-valued quantity, regardless if it is real, complex, quaternions, etc. Indeed, real ($\mathbb{R}$), complex ($\mathbb{C}$), and quaternion algebras ($\mathbb{H}$), commonly used in adaptive filtering and optimization literature~\cite{MandicBook,Took_QLMS,Mandic2011,Neto_WLQLMS_2011,Jahanchahi2012,Mengdi2014,Dongpo2016}, are subsets (subalgebras) of the GA created by general multivectors like~\eqref{def:multivector}. Thus, to support the generalization of standard AFs achieved by GAAFs (Section~\ref{cap:GAAFs_standard}), this section shows how those subalgebras of interest can be retrieved from $\mathcal{G}(\mathbb{R}^3)$, the complete GA of $\mathbb{R}^3$, via isomorphisms\footnote{For simulation purposes, this paper focuses on the case $n=3$, i.e., the subalgebras of $\mathcal{G}(\mathbb{R}^3)$. However, the GAAFs can work with any subalgebra of $\mathcal{G}(\mathbb{R}^n), \forall n\in \mathbb{Z}\setminus\{0\}$ (Section~\ref{cap:GAAFs_standard}).}.

\begin{figure}[!t]
	\centering
	\includegraphics[width=0.39\columnwidth, angle=-90]{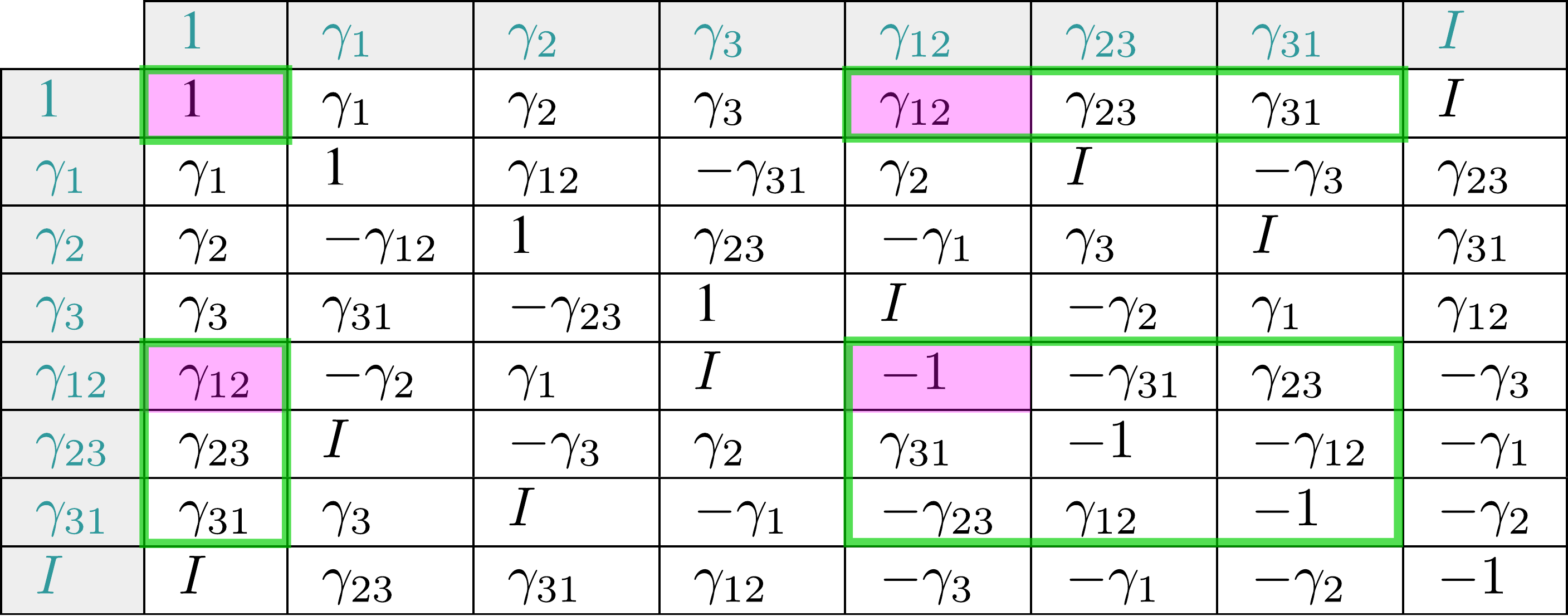}
	\caption{Multiplication table of $\mathcal{G}(\mathbb{R}^3)$ via the geometric product. The elements highlighted in
		magenta constitute the multiplication table of $\mathcal{G}^{+}(\mathbb{R}^2)$ (even subalgebra isomorphic to complex numbers). Similarly, the elements enclosed in the green rectangles form the multiplication table of $\mathcal{G}^{+}(\mathbb{R}^3)$ (even subalgebra isomorphic to quaternions).}
	\label{fig:multiplication_table}
\end{figure}

The \textit{Complete Geometric Algebra of} $\mathbb{R}^3$ is obtained by multiplying the elements of~\eqref{eq:orthonormal_basis_R3} via the geometric product. As depicted in \figurename~\ref{fig:multiplication_table}, this results in several terms which linearly combined represent all the possible multivectors in $\mathcal{G}(\mathbb{R}^3)$. Also, note that since $\mathbb{R}^2 \subset \mathbb{R}^3 \Rightarrow \mathcal{G}(\mathbb{R}^2) \subset  \mathcal{G}(\mathbb{R}^3)$, the multiplication table for $\mathcal{G}(\mathbb{R}^2)$ can be recovered from \figurename~\ref{fig:multiplication_table}.

The \textit{even subalgebras}~\cite{Hitzer_Introduction,hestenes1999newfoundations,2016introductionVaz} of $\mathcal{G}(\mathbb{R}^3)$ and $\mathcal{G}(\mathbb{R}^2)$, i.e., those whose elements have only even grades ($g = 0,2,4,\cdots$ in~\eqref{def:multivector}) and denoted $\mathcal{G}^{+}(\mathbb{R}^3)$ and $\mathcal{G}^{+}(\mathbb{R}^2)$, are of special interest. One can show that $\mathcal{G}^{+}(\mathbb{R}^2)$, with basis $\{1,\gamma_{12}\}$ (also $\{1,\gamma_{23}\}$ and $\{1,\gamma_{31}\}$), is \textit{isomorphic to the complex numbers}~\cite{Hitzer_Introduction}. This is established by identifying the imaginary unit $j$ with the bivector $\gamma_{12}$, $j = \gamma_{12} = \gamma_1\gamma_2 = \gamma_1\wedge\gamma_2$. From \figurename~\ref{fig:multiplication_table} it is known that $(\gamma_{12})^2 = -1$. Then, $j^2 = (\gamma_{12})^2 = -1$, demonstrating the isomorphism. Similarly, $\mathcal{G}^{+}(\mathbb{R}^3)$ with basis $\{1,\gamma_{12},\gamma_{23},\gamma_{31}\}$ is shown to be \textit{isomorphic to quaternion algebra} via the adoption of the following correspondences: $i \leftrightarrow -\gamma_{12}$, $j \leftrightarrow -\gamma_{23}$, $k \leftrightarrow -\gamma_{31}$, where $\{i,j,k\}$ are the three imaginary unities of a quaternion. The minus signs are necessary to make the product of two bivectors equal to the third one and not minus the third, e.g., $(-\gamma_{12})(-\gamma_{23}) = \gamma_{13} = -\gamma_{31}$, just like quaternions, i.e. $ij = k$, $jk = i$, and $ki = j$~\cite{Hitzer_Introduction,quaternionReport,Girard2007}. Additionally, note that $\mathcal{G}^{+}(\mathbb{R}^2) \subset \mathcal{G}^{+}(\mathbb{R}^3) \subset \mathcal{G}(\mathbb{R}^3)$.

It follows that the dimension of the complete GA of $\mathcal{G}(\mathcal{V})$ can be obtained from its subalgebras, i.e., $\dim\{\mathcal{G}(\mathcal{V})\} = \sum\limits^{n}_{g=0} \dim\{\mathcal{G}^g(\mathcal{V})\} = \sum\limits^{n}_{g=0} \binom n{g} = 2^n$. 

\subsection{Performing Rotations with Multivectors (Rotors)}
\label{ssec:performing_rotation}
The even subalgebra $\mathcal{G}^{+}(\mathbb{R}^n)$ is also known as the \emph{algebra of rotors}, i.e., its elements are a special type of multivector (called rotors) able to \textit{apply rotations} to vectors in $\mathbb{R}^n$~\cite{hestenes1987GAcalculus,hestenes1999newfoundations}. Given a vector $x \in \mathbb{R}^n$, it can be rotated by applying the rotation operator $r(\cdot)\t~r$
\vspace*{-5mm}
\eq 
x \rightarrow \underbrace{rx\t~r}_{rotated},
\label{eq:rotation_operator}
\vspace*{-2mm}
\qe     
where $r \in \mathcal{G}^{+}(\mathbb{R}^n)$ is a rotor, $\t~r$ is its reverse, and $r\t~r = 1$, i.e., $r$ is a unit rotor. Note that the unity constraint is necessary to avoid the rotation operator to scale the vector $x$, i.e., to avoid changing its norm. 

\begin{figure}[!t]
	\centering
	\subfloat[]{\includegraphics[trim=0cm 0.3cm 0cm 0cm, clip=true,width=0.12\textwidth]{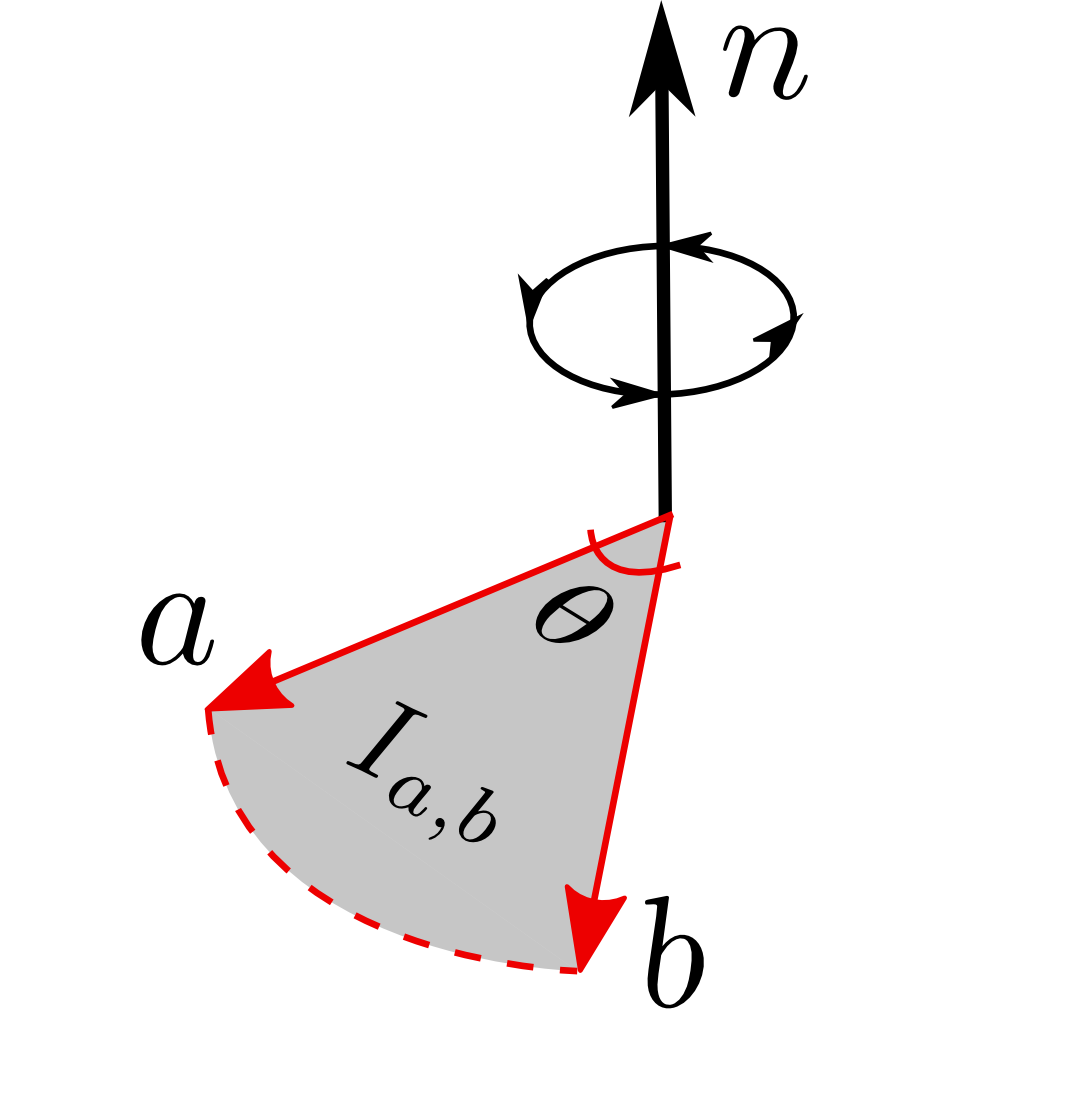}}
	~\hspace*{8mm}
	\subfloat[]{\includegraphics[trim=0cm 0.3cm 0cm 0cm, clip=true,width=0.12\textwidth]{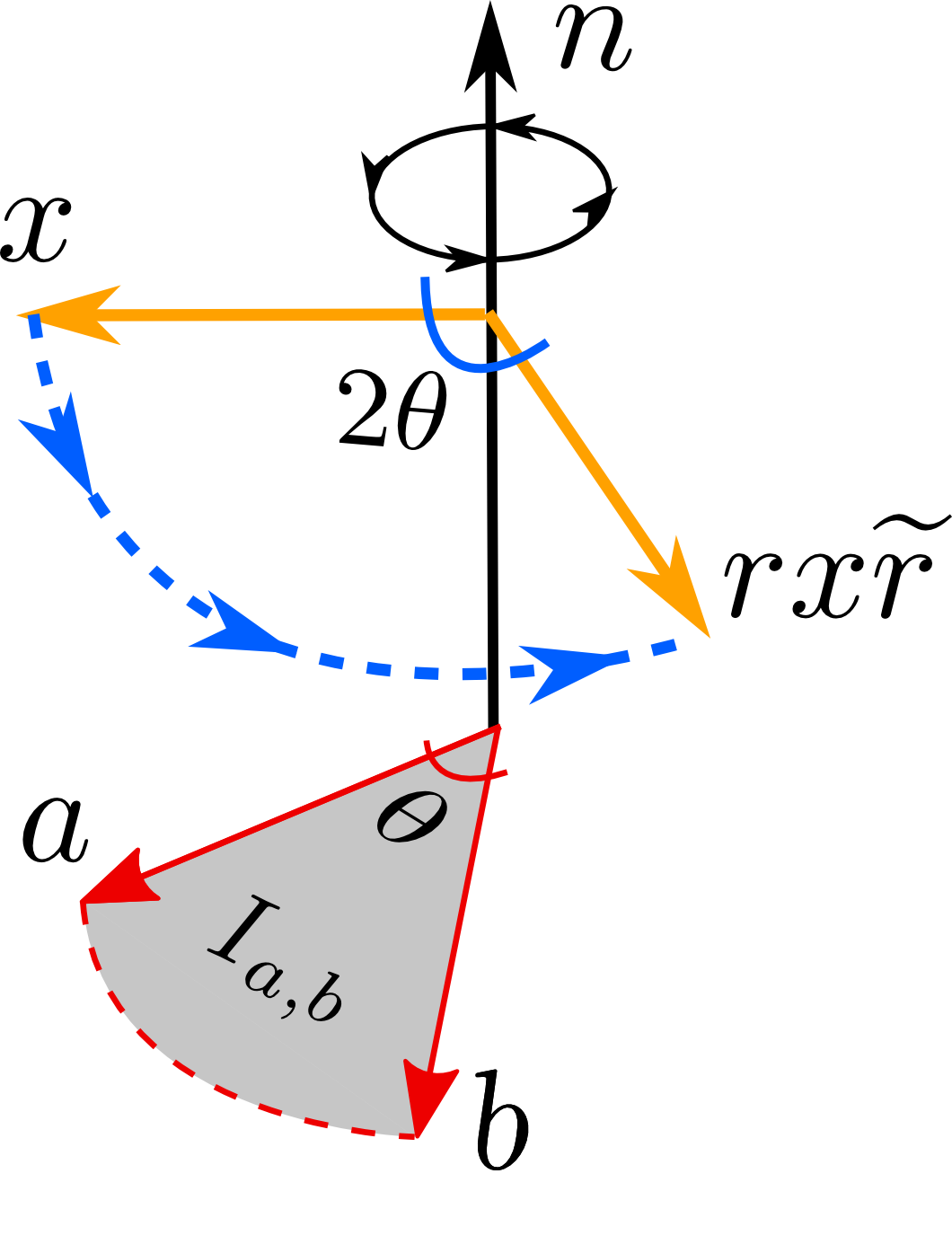}}
	\vspace*{-1mm}
	\caption[Rotation operator]{(a) A rotor can be generated from the geometric multiplication of two unit vectors in $\mathbb{R}^n$. (b) Applying the rotation operator: the vector $x$ is rotated by an angle of $2\theta$ about the normal $n$ of the oriented surface $I_{a,b}$ (a hyperplane). In two dimensions (2D), $r$ is a complex number, $\t~r$ its conjugate, and they rotate $x$ about the normal of the complex plane.}
	\label{fig:rotation_operator}
\end{figure}

A rotor $r \in \mathcal{G}^{+}(\mathbb{R}^n)$ can be generated from the geometric multiplication of two unit vectors in $\mathbb{R}^n$. Given $\{a,b\} \in \mathbb{R}^n$, $|a| = |b| = 1$, with an angle $\theta$ between them, and using~\Crefrange{def:inner_product_vectors}{def:geometric_product_vectors}, the \emph{exponential form of a rotor} is~\cite[p. 107]{hestenes1987GAcalculus} 
\vspace*{-2mm}
\eq 
\eqarray
r = ab = a\cdot b + a\wedge b \hspace*{-2mm}&= |a||b|cos\theta + I_{a,b}|a||b|sin\theta\\	
&= cos\theta + I_{a,b}sin\theta = e^{I_{a,b}\theta}.
\qearray 
\label{eq:rotor_definition}
\vspace*{-1mm}
\qe %
As shown in \figurename~\ref{fig:rotation_operator}, $x$ is rotated by an angle of $2\theta$ about the normal of the oriented surface $I_{a,b}$ (rotation axis)~\cite{hestenes1999newfoundations}. Note that both quantities ($\theta$ and $I_{a,b}$) define the rotor in~\eqref{eq:rotor_definition}. The geometric transformation enabled by rotors was applied in~\cite{SPL,WACV} to devise AFs that estimate the relative rotation between 3D Point Clouds (PCDs).   

\section{Linear Estimation in GA}\label{cap:linear_estimation_GA}

This section introduces an instantaneous quadratic cost function with multivector entries. This is key to expand the estimation capabilities of AFs, generating the GAAFs further ahead in Section~\ref{cap:GAAFs_standard}. Moreover, it is shown that the standard LA counterpart can be recovered as a special case.   

\subsection{Definitions}
\label{sec:useful_defs_LinEst_in_GA}

The \emph{scalar product} between two multivectors is $A * B \triangleq \langle AB \rangle$, i.e., it is the scalar part ($0$-grade) of the geometric multiplication between $A$ and $B$ (for the special case of vectors, $a * b = \langle ab \rangle = a \cdot b$). Its commutativity originates the \emph{cyclic reordering property}~\cite{hestenes1987GAcalculus} $A * B = \langle AB \rangle = \langle BA \rangle = B * A \Rightarrow  \langle AB \cdots C \rangle = \langle B \cdots CA \rangle$. 

An \emph{array of multivectors} is a collection (row or column) of general multivectors. Given $M$ multivectors $\{U_1, U_2, \cdots, U_M\}$ in $\mathcal{G}(\mathbb{R}^3)$, the $M\times 1$ array collects them as follows
{
	\begin{align}
	u &{=} \begin{bmatrix}
	U_1 \\
	\vdots \\
	U_M
	\end{bmatrix}{=} 
	\begin{bmatrix}
	u{\scriptstyle(1,0)} {+} u{\scriptstyle(1,1)}\gamma_1 {+} \cdots {+} u{\scriptstyle(1,6)}\gamma_{31} {+} u{\scriptstyle(1,7)}I \\
	\vdots \\
	u{\scriptstyle(M,0)} {+} u{\scriptstyle(M,1)}\gamma_1 {+} \cdots {+} u{\scriptstyle(M,6)}\gamma_{31} {+} u{\scriptstyle(M,7)}I \\
	\end{bmatrix}. 
	\label{eq:array_MV_example_1}
	\end{align}
}%
The array is denoted using lower case letters, the same as scalars and vectors (1-vectors). However, the meaning of the symbol will be evident from the context. Additionally, this work adopts the notion of \textit{matrix of multivectors} as well: it is a matrix whose elements are multivector-valued.  

Given two $M\times 1$ arrays of multivectors, $u$ and $w$, the \emph{array product} between them is defined as
\vspace*{-2mm} 
{\small 
	\begin{align}
	u^{T}w &= 
	U_1W_1 +
	U_2W_2 + 
	\cdots
	U_MW_M = \sum\limits_{j=1}^{M} U_jW_j,
	\label{eq:array_product_1}
	\end{align}
}%
which results in a general multivector. The underlying product in each of the terms $U_jW_j$, $j = \{1,\cdots, M\}$, is the geometric product. Observe that due to the noncommutativity of the geometric product, $u^{T}w \neq w^{T}u$ in general. 

Similarly, multiplications involving matrices of multivectors follow the general rules of matrix algebra, however using the geometric product as the underlying product, just like in~\eqref{eq:array_product_1}.   

The \emph{reverse transpose array} is the extension of the reverse operation of multivectors to include arrays of multivectors. Given the array $u$ in~\eqref{eq:array_MV_example_1}, its reverse version, denoted by $\rt{}$, is
\vspace*{-3mm}
{\normalsize 
	\begin{align}
	\rt{u} &= \begin{bmatrix}
	\t~U_1\hspace*{1mm}\t~U_2\hspace*{1mm}\cdots\hspace*{1mm}\t~U_M 
	\end{bmatrix}. 
	\label{eq:reverse_transposed_array}
	\end{align}
}%
Note the similarity with the Hermitian conjugate, its counterpart in LA. From~\eqref{eq:array_product_1} and~\eqref{eq:reverse_transposed_array} it follows that
\vspace*{-2mm}
{\small 
	\eq 
	\rt{u}w = \sum\limits_{j=1}^{M} \t~U_jW_j,
	\label{eq:array_product}
	\vspace*{-2mm}
	\qe 
}%
which results in a general multivector.

The array product between $\rt u$ and $u$ is represented by the notation $\norm{u}^2 \triangleq \rt u u$. Note the same notation is employed in LA to represent the squared norm of a vector in $\mathbb{R}^n$. However, from~\eqref{eq:array_product} it is known that $\norm{u}^2 = \rt u u$ \emph{is a general multivector}, i.e., \emph{it is not a pure scalar value} which in LA provides a measure of distance. In GA, the distance metric is given by the \emph{magnitude of a multivector}, defined in terms of the scalar product, e.g., $|A| \triangleq \sqrt{A * \widetilde{A}} = \sqrt{\sum_{g} | \langle A \rangle_g|^2}\text{,}$ which is indeed a scalar value. Thus, for an array $u$ and a multivector $U$, 
\vspace*{-2mm}
\eq
\eqarray 
&\norm{u}^2 = \rt u u : \text{is a multivector}\\
&|U|^2 = U {*} \widetilde{U} = \widetilde{U} {*} U = \sum_{g} | \langle U \rangle_g|^2 : \text{is a scalar}. 
\qearray 
\label{eq:norm_of_array}
\vspace*{-2mm}
\qe

The product between a multivector $U$ and an array $w$, namely $Uw$, is defined as the geometric multiplication of $U$ by each entry of $w$ (a procedure similar to multiplying a scalar by a vector in LA). Due to the noncommutativity of the geometric product, $Uw \neq wU$ in general. 

\subsection{General Cost Function in GA}
\label{sec:CF_in_GA}

Following the guidelines in~\cite[p.64 and p.121]{hestenes1987GAcalculus}, one can formulate a minimization problem in GA by defining the following cost function
{\small
	\vspace*{-3mm}
	\eq
	J(D,A_k,X,B_k) = \left\lvert D - \sum\limits_{k=1}^{M} A_k X B_k \right\rvert^2,
	\label{eq:general_CF_GA} 
	\vspace*{-2mm}
	\qe
}
where $D, A_k, X, B_k$ are general multivectors. The term $\sum\limits_{k=1}^{M} A_k X B_k$ (the addition of $M$ multiplications) represents the canonical form of a multilinear transformation applied to the multivector $X$~\cite[p.64 and p.121]{hestenes1987GAcalculus}. The goal is to optimize $\{A_k, B_k\}$ in order to minimize~\eqref{eq:general_CF_GA}.

Special cases of~\eqref{eq:general_CF_GA} are obtained depending on the values chosen for $D, A_k, X, B_k$ and $M$. For instance, making $M = 1$, $D = y$ (1-vector), $X = x$ (1-vector), $A_k = r$ (rotor), and $B_k = \t~r$ yields $J(r) = \left\lvert y - rx\t~r\right\rvert^2$, the instantaneous version of the cost function minimized in~\cite{SPL,WACV} (subject to $r\t~r = \t~rr = 1$) to estimate the relative rotation between 3D PCDs. 

In this paper it is studied the case in which $X = \t~{U}_k$, $A_k = 1$, $B_k = W_k$ (general multivectors), so that   
\vspace*{-1mm}
{\small   
	\eq
	J(w) = \left\lvert D - \sum\limits_{k=1}^{M} \t~U_k W_k \right\rvert^2 = \left\lvert D - \rt{u}w\right\rvert^2 = \left\lvert E \right\rvert^2, 
	\label{eq:AF_cost_function_1}
	\vspace*{-2mm}
	\qe 
}
where $M$ now also represents the system order (the number of taps in the filter), $E = D - \rt{u}w$ is the estimation error, and the definition of array product~\eqref{eq:array_product} was employed to make $\sum_{k=1}^{M} \t~U_k W_k = \rt{u}w$. Notice that \eqref{eq:AF_cost_function_1} is the GA counterpart of the standard LA cost function~\cite[p.477]{Sayed08}. And similarly to its LA counterpart, $D$ (multivector) is estimated as a linear combination of the multivector entries of $u$.

\section{Geometric-Algebra Adaptive Filters}\label{cap:GAAFs_standard}


In this section, the GAAFs are motivated by deriving the GA-LMS to minimize the cost function~\eqref{eq:AF_cost_function_1} in an adaptive manner. This is done
by writing~\eqref{eq:AF_cost_function_1} at instant $i$, yielding an instantaneous cost function $J(i)$~\cite{Diniz2013,HaykinWidrow}, as shown in the sequel.      


At each iteration $i$, two multivector-valued signals $D(i)$ and $U(i)$
are observed. Samples of $U(i)$ are collected into an $M\times 1$ array $\rt{u}_i = \begin{bmatrix} \t~U_1(i)\hspace*{1mm}\t~U_2(i)\hspace*{1mm}\cdots\hspace*{1mm}\t~U_M(i) \end{bmatrix}$.
The linear combination (array product) $\hat{D}(i) = \rt{u}_i w_{i-1}$ is used as an estimate of $D(i)$, generating the \textit{a priori} output estimation error 
\vspace*{-2mm}
\eq
E(i) = D(i) - \hat{D}(i) = D(i) - \rt{u}_iw_{i-1}.
\label{eq:estimation_error}
\vspace*{-2mm}
\qe 
This way, from~\eqref{eq:AF_cost_function_1} and~\eqref{eq:estimation_error}
\vspace*{-2mm}
\eq
J(w_{i-1}) = \left\lvert D(i) - \rt{u}_iw_{i-1}\right\rvert^2 = \left\lvert E(i) \right\rvert^2. 
\label{eq:AF_cost_function_1_again}
\vspace*{-2mm}
\qe 


Given $E(i)$ and $w_{i-1}$, the GAAFs update the estimate for the array of multivectors $w$ via a recursive rule of the form
\vspace*{-2mm}
\begin{equation}
w_i = w_{i-1} + \mu h \text{,}
\label{eq:GAAF_update}
\vspace*{-2mm}
\end{equation}  
where $\mu$ (a scalar) is the AF step size (learning rate), and $h$ is an array of multivectors related to the estimation error $e(i)$. This work adopts an instantaneous steepest-descent rule~\cite{Sayed08,Diniz2013,HaykinWidrow}, in which the AF is designed to follow the opposite direction of the instantaneous gradient of~\eqref{eq:AF_cost_function_1_again}, namely $\partial_w J(w_{i-1})$. This way, $h \triangleq -B\partial_w J(w_{i-1})$, yielding the general form of a GA-based adaptive rule
\vspace*{-2mm}
\eq
w_{i} = w_{i-1} - \mu B \left[\partial_w J(w_{i-1})\right],
\label{eq:AF_general_form}
\vspace*{-2mm}
\qe  
in which $B$ is a matrix with multivector entries. Choosing $B$ appropriately is a requirement to define the type of adaptive algorithm~\cite{Sayed08}.   


\subsection{A Note on Multivector Derivative}
\label{ssec:multivector_derivative}
Equation~\eqref{eq:AF_general_form} requires the calculation of a \emph{multivector derivative}. In GA, the differential operator $\partial_w$ has the algebraic properties of a multivector in $\mathcal{G}(\mathbb{R}^n)$~\cite{Hitzer_MultivectorDiffCalc}. Put differently, the gradient $\partial_w J(w_{i-1})$ can be interpreted as the \emph{geometric multiplication} between the multivector-valued quantities $\partial_w$ and $J(w_{i-1})$. 

This way, it follows from~\eqref{def:multivector} that $\partial_w$ can be decomposed into its basis blades. In fact, it is known that any multivector $A \in \mathcal{G}(\mathbb{R}^n)$ can be decomposed into blades via
\begin{equation}
A = \sum\limits_{k=0}^{2^n - 1} \basis_k (\basis^k * A) = \sum\limits_{k=0}^{2^n - 1} \basis_k \langle \basis^k A \rangle = \sum\limits_{k=0}^{2^n - 1} \basis_k A_k,
\label{eq:decomposition_grades}
\end{equation}
in which $A_k$ is scalar-valued, and $\{\basis^k\}$ and $\{\basis_k\}$, $k = 0,\cdots,2^n{-}1$ are the so-called \textit{reciprocal (dual) bases} of $\mathcal{G}(\mathbb{R}^n)$\footnote{Two symbols are used to refer to GA basis, each with a different purpose. The symbol $\gamma_{ij}$ (with two indexes) is adopted when dealing with specific algebras, e.g., $\mathcal{G}(\mathbb{R}^2)$ or $\mathcal{G}(\mathbb{R}^3)$. For general algebras $\mathcal{G}(\mathbb{R}^n)$ and analytical derivation, the symbol $\basis_{k}$ (with only one index) is more appropriate.}. The concept of reciprocal (dual) bases is a useful analytical tool in GA to convert from nonorthogonal to orthogonal vectors and vice versa -- it simplifies the analytical procedure ahead since mutually orthogonal elements cancel out. Details are provided in~\cite[Section 1-3]{hestenes1987GAcalculus} and~\cite{Hitzer_Introduction}. For the purpose of this paper, it suffices to know that the following relation holds for dual bases: $\basis^k * \basis_j = \langle \basis^k\basis_j \rangle = \delta_{k,j}$, where $\delta_{k,j} = 1$ for $k=j$ (\textit{Kronecker delta}). It is easy to show that the basis for $\mathcal{G}(\mathbb{R}^n)$, $\{\basis_k\}, \text{for }k=0,\cdots\hspace*{-0.5mm},2^n{-}1$ and its reversed version $\{\t~{\basis_k}\}$ comply with the relation above (i.e., they are dual bases). Therefore, they are utilized from now on to decompose multivectors into blades. In particular, applying~\eqref{eq:decomposition_grades} to $\partial_w$ results in
\vspace*{-2mm}
\eq
\partial_w \triangleq \sum_{\ell=0}^{2^n{-}1} \basis_\ell \langle \widetilde{\basis}_\ell \partial_{w} \rangle = \sum_{\ell=0}^{2^n{-}1} \basis_\ell \partial_{w,\ell}, 
\label{eq:diff_operator}
\vspace*{-2mm}
\qe  
where each term $\partial_{w,\ell}$ in the sum is the usual derivative from standard calculus, which affects \emph{only} blade $\ell$. Form~\eqref{eq:diff_operator} provides some analytical advantages (see~\cite{Hitzer_split}) and is employed next to calculate the gradient $\partial_w J(w_{i-1})$.

\subsection{Calculating the Multivector-valued Gradient}
\label{ssec:multivector-valued_gradient}

Noticing that~\eqref{eq:AF_cost_function_1_again} can be written in terms of the scalar product
\vspace*{-2mm}
\eq
J(w_{i-1}) = \left\lvert E(i) \right\rvert^2 = E(i)*\widetilde{E}(i),
\label{eq:AF_cost_function_scalar_product}
\vspace*{-2mm}
\qe 
one can decompose it \emph{in terms of its blades} via~\eqref{eq:decomposition_grades},
\vspace*{-2mm}
\eq
\eqarray
J(w_{i-1}) \hspace*{-2mm}&= \Big(\sum_{\mathcal{p}=0}^{2^n{-}1} \basis_\mathcal{p} E_\mathcal{p}\Big) * \Big(\sum_{\mathcal{p}=0}^{2^n{-}1} E_\mathcal{p} \widetilde{\basis}_\mathcal{p}\Big) \\ &= \sum_{\mathcal{p}=0}^{2^n{-}1} E^2_\mathcal{p} (\basis_\mathcal{p} {*} \t~\basis_\mathcal{p}) {=} \sum_{\mathcal{p}=0}^{2^n{-}1} E^2_\mathcal{p} \langle\basis_\mathcal{p}\t~\basis_\mathcal{p}\rangle {=} \sum_{\mathcal{p}=0}^{2^n{-}1} E^2_\mathcal{p}.
\qearray
\label{eq:cost_function_freqtrack_grades} 
\qe

The gradient $\partial_w J(w_{i-1})$ is obtained by multiplying~\eqref{eq:diff_operator} and~\eqref{eq:cost_function_freqtrack_grades} 
\vspace*{-2mm}
\eq
\partial_w J(w_{i-1}) {=} \Bigg(\hspace{-1mm}\sum_{\ell=0}^{2^n{-}1} \basis_\ell \partial_{w,\ell}\Bigg) \Bigg(\sum_{\mathcal{p}=0}^{2^n{-}1} E^2_\mathcal{p}\Bigg) {=} \hspace{-2mm}\sum_{\mathcal{p},d=0}^{2^n{-}1}\hspace{-2mm}\basis_\ell \partial_{w,\ell} E^2_\mathcal{p}.
\label{eq:grad_1}
\vspace*{-2mm}
\qe
From~\eqref{eq:estimation_error} one can notice that $E_\mathcal{p} \triangleq D_\mathcal{p} - \hat{D}_\mathcal{p}$, where $\hat{D}_\mathcal{p} = \langle \hat{D}(i)\rangle_\mathcal{p} = \langle \rt{u}_iw_{i-1}\rangle_\mathcal{p}$ (similar for $D_\mathcal{p}$). Thus
\vspace*{-2mm} 
\eq
\eqarray
\partial_{w,\ell} E^2_\mathcal{p} {=} 2 E_\mathcal{p} (\partial_{w,\ell} E_\mathcal{p}) {=} 2 E_\mathcal{p} (\partial_{w,\ell} (D_\mathcal{p} {-} \hat{D}_\mathcal{p} )) {=} {-}2 E_\mathcal{p} (\partial_{w,\ell}\hat{D}_\mathcal{p}),
\qearray
\label{eq:grad_2}
\vspace*{-1mm}
\qe
where $\partial_{w,\ell}D_\mathcal{p} = 0$ since $D_\mathcal{p}$ does not depend on the weight array $w$.

Plugging \eqref{eq:grad_2} into \eqref{eq:grad_1} results in
\vspace*{-2mm}
\eq
\partial_w J(w_{i-1}) = -2\hspace{-2mm}\sum_{\mathcal{p},\ell=0}^{2^n{-}1} \basis_\ell E_\mathcal{p} (\partial_{w,\ell}\hat{D}_\mathcal{p}).
\label{eq:grad_3}
\vspace*{-2mm} 
\qe

The term $\hat{D}_\mathcal{p}$ is obtained by decomposing $\hat{D}(i)$ into its blades   
\vspace*{-2mm}
\eq
\hat{D}(i) = \rt u_{i}\iprod w_{i-1} = \sum_{\mathcal{p}=0}^{2^n{-}1} \basis_\mathcal{p} \langle \widetilde{\basis}_\mathcal{p} (\rt u_{i}\iprod w_{i-1}) \rangle,
\label{eq:def_d_hat}
\vspace*{-2mm}
\qe 
which requires to perform the decomposition of $\rt u_{i}$ and $w_{i-1}$ (arrays of multivectors). Indeed, arrays of multivectors can be decomposed into blades as well. For instance 
\begin{align}
\begin{bmatrix}
1 + \gamma_{12} \\
2 + \gamma_{1} + \gamma_{3} \\
\end{bmatrix}{=} 
\begin{bmatrix}
1 \\
2 \\
\end{bmatrix}{+} 
\begin{bmatrix}
0 \\
1 \\
\end{bmatrix}{\gamma_{1}}{+} 
\begin{bmatrix}
0 \\
1 \\
\end{bmatrix}{\gamma_{3}}{+} 
\begin{bmatrix}
1 \\
0 \\
\end{bmatrix}{\gamma_{12}}. 
\label{eq:decomposing_array}
\end{align}
Thus, employing~\eqref{eq:decomposition_grades} once again, $u_{i}$ and $w_{i-1}$ can be written in terms of their $2^n$ blades 
\vspace*{-2mm}
\eq
\rt u_{i} = \sum_{\mathcal{q}=0}^{2^n{-}1} u^{T}_{i,\mathcal{q}} \t~\basis_\mathcal{q} \text{\hspace*{2mm}and\hspace*{2mm}} w_{i-1} = \sum_{\mathcal{t}=0}^{2^n{-}1} \basis_\mathcal{t} w_{i-1,\mathcal{t}},  
\label{eq:u_grades_and_w_grades} 
\qe
where $u^{T}_{i,\mathcal{q}}$ and $w_{i-1,\mathcal{t}}$ are respectively $1\times M$ and $M \times 1$ arrays with real entries. Plugging~\eqref{eq:u_grades_and_w_grades} back into~\eqref{eq:def_d_hat}\footnote{From now on, the iteration subscripts $i$ and $i-1$ are omitted from $u_{i,\mathcal{q}}$ and $w_{i-1,\mathcal{t}}$ for clarity purposes.}
\eq
\eqarray
\hat{D}(i) = \rt u_{i}\iprod w_{i-1} &\hspace{-3mm}= \sum_{\mathcal{p}=0}^{2^n{-}1} \basis_\mathcal{p} \langle \widetilde{\basis}_\mathcal{p} ( \sum_{\mathcal{q}=0}^{2^n{-}1} u^{T}_\mathcal{q} \t~\basis_\mathcal{q}  \iprod \sum_{\mathcal{t}=0}^{2^n{-}1} \basis_\mathcal{t} w_\mathcal{t})\rangle \vspace{2mm} \\
&\hspace{-3mm}= \sum_{\mathcal{p}=0}^{2^n{-}1} \basis_\mathcal{p} \sum_{\mathcal{q},\mathcal{t}=0}^{2^n{-}1} \langle \widetilde{\basis}_\mathcal{p} (u^{T}_\mathcal{q} \t~\basis_\mathcal{q} \iprod \basis_\mathcal{t} w_\mathcal{t}) \rangle \vspace{2mm}  \\
&\hspace{-3mm}= \sum_{\mathcal{p}=0}^{2^n{-}1} \basis_\mathcal{p} \sum_{\mathcal{q},\mathcal{t}=0}^{2^n{-}1} \langle \widetilde{\basis}_\mathcal{p}\t~\basis_\mathcal{q}\basis_\mathcal{t}\rangle (u^{T}_\mathcal{q} \cdot w_\mathcal{t})\vspace{2mm} \\
&\hspace{-3mm}= \sum_{\mathcal{p}=0}^{2^n{-}1} \basis_\mathcal{p} \hat{D}_\mathcal{p}, \text{ in which}
\qearray
\vspace*{-1mm}
\qe 
\eq
\hat{D}_\mathcal{p} = \sum_{\mathcal{q},\mathcal{t}=0}^{2^n{-}1} \langle \widetilde{\basis}_\mathcal{p}\t~\basis_\mathcal{q}\basis_\mathcal{t}\rangle (u^{T}_\mathcal{q} \cdot w_\mathcal{t}), \hspace*{1mm}\mathcal{p}=0,\cdots,2^n{-}1
\label{eq:da_coefficients}
\vspace*{-2mm}
\qe
is the expression of $\hat{D}_\mathcal{p}$ as a function of the blades of $\rt u_{i}\iprod w_{i-1}$.

From~\eqref{eq:da_coefficients}, the term $\partial_{w,\ell}\hat{D}_\mathcal{p}$ of~\eqref{eq:grad_3} becomes   
\vspace*{-2mm} 
\eq
\eqarray
\partial_{w,\ell}\hat{D}_\mathcal{p} &\hspace{-2mm}{=} \partial_{w,\ell} \left[ \sum_{\mathcal{q},\mathcal{t}=0}^{2^n{-}1} \langle \widetilde{\basis}_\mathcal{p}\t~\basis_\mathcal{q}\basis_\mathcal{t}\rangle (u^{T}_\mathcal{q} \cdot w_\mathcal{t}) \right]\\
&\hspace{-2mm}{=} \sum_{\mathcal{q},\mathcal{t}=0}^{2^n{-}1} \langle \widetilde{\basis}_\mathcal{p}\t~\basis_\mathcal{q}\basis_\mathcal{t}\rangle \partial_{w,\ell} (u^{T}_\mathcal{q} \cdot w_\mathcal{t}).
\qearray
\label{eq:partial_w_aux_1}
\vspace*{-2mm} 
\qe

It is important to notice that $\partial_{w,\ell} (u^{T}_\mathcal{q} \cdot w_\mathcal{t})$ will be different from zero only when $\ell = \mathcal{t}$, i.e., when $\partial_{w,\ell}$ and $w_\mathcal{t}$ are of same blade. This is the case since $\partial_{w,\ell}$ is the partial derivative of $w$ with respect to the blade $\ell$ \emph{only}. Therefore, if $\ell \neq \mathcal{t}$ then the partial derivation yields zero, i.e., $\partial_{w,\ell} w_\mathcal{t} = 0 \Rightarrow \partial_{w,\ell} (u^{T}_\mathcal{q} \cdot w_\mathcal{t}) = 0$ (note that $u^{T}_\mathcal{q}$ does not depend on $w$). Thus, adopting the \textit{Kronecker delta function}~\cite{Hitzer_split}, $\partial_{w,\ell} (u^{T}_\mathcal{q} \cdot w_\mathcal{t}) = \delta_{\mathcal{t},\ell} u^{T}_\mathcal{q}$, and~\eqref{eq:partial_w_aux_1} becomes
\vspace*{-2mm} 
\eq 
\partial_{w,\ell}\hat{D}_\mathcal{p} = \sum_{\mathcal{q},\mathcal{t}=0}^{2^n{-}1} \langle \widetilde{\basis}_\mathcal{p}\t~\basis_\mathcal{q}\basis_\mathcal{t}\rangle\delta_{\mathcal{t},\ell} u^{T}_\mathcal{q}.
\label{eq:differential_da}
\vspace*{-2mm} 
\qe

Finally, substituting~\eqref{eq:differential_da} into~\eqref{eq:grad_3}, the stochastic gradient is obtained
\vspace*{-1mm}
{\small 
	\eq
	\eqarray
	\partial_w J(w_{i-1}) &\hspace{-2mm}{=} -2\sum_{\mathcal{p},\ell=0}^{2^n{-}1} \basis_\ell E_\mathcal{p} \sum_{\mathcal{q},\mathcal{t}=0}^{2^n{-}1} \langle \widetilde{\basis}_\mathcal{p}\t~\basis_\mathcal{q}\basis_\mathcal{t}\rangle\delta_{\mathcal{t},\ell} u^{T}_\mathcal{q} \vspace{2mm} \\
	&\hspace{-2mm}{=} -2\sum_{\mathcal{p},\ell=0}^{2^n{-}1} E_\mathcal{p} \sum_{\mathcal{q}=0}^{2^n{-}1} \basis_\ell  \langle \widetilde{\basis}_\mathcal{p}\t~\basis_\mathcal{q}\basis_\ell\rangle u^{T}_\mathcal{q} \vspace{2mm} \\
	&\hspace{-2mm}{=} -2\sum_{\mathcal{p},\ell=0}^{2^n{-}1} E_\mathcal{p} \basis_\ell  \langle \widetilde{\basis}_\mathcal{p}\rt u_{i}\basis_\ell\rangle \vspace{2mm} \\
	&\hspace{-2mm}{=} -2\sum_{\mathcal{p},\ell=0}^{2^n{-}1} E_\mathcal{p} \basis_\ell  \langle \widetilde{\basis}_\ell u_{i}\basis_\mathcal{p}\rangle \vspace{2mm} \\
	&\hspace{-2mm}{=} -2\sum_{\mathcal{p}=0}^{2^n{-}1} E_\mathcal{p} u_{i}\basis_\mathcal{p} = \boxed{-2 u_{i}E(i)}. 
	\qearray
	\label{eq:grad_4}
	\vspace*{-1.5mm} 
	\qe 
}

In the AF literature, setting $B$ equal to the identity matrix in~\eqref{eq:AF_general_form} results in the stochastic-gradient update rule~\cite{Sayed08}. This is adopted here as well in order to devise the GA version of the LMS filter -- however, GA allows for selecting $B$ with multivector entries, opening up the possibility to generate other types of GA-based adaptive filters. Substituting~\eqref{eq:grad_4} into~\eqref{eq:AF_general_form} and setting $B$ equal to the identity matrix yields the GA-LMS update rule
\vspace*{-2mm} 
\eq
\boxed{w_{i} = w_{i-1} + \mu u_{i} e(i)},
\label{eq:GA-LMS_standard_updateRule} 
\vspace*{-1.5mm} 
\qe
where the $2$ in~\eqref{eq:grad_4} was absorbed by the scalar step size $\mu$.

Note that the GA-LMS~\eqref{eq:GA-LMS_standard_updateRule} has the same format of standard LMS AFs~\cite{Diniz2013}, namely the \emph{real-valued} LMS ($u$ and $w$ have real-valued entries) and the \emph{complex-valued} LMS ($u$ and $w$ have complex-valued entries). However, this work puts no constraints on the entries of the arrays $u$ and $w$ -- they can be any kind of multivector. This way, the update rule~\eqref{eq:GA-LMS_standard_updateRule} is valid for any $u$ and $w$ whose \emph{entries are general multivectors} in $\mathcal{G}(\mathbb{R}^n)$. In other words, \eqref{eq:GA-LMS_standard_updateRule} generalizes the standard LMS AF for several types of $u$ and $w$ entries: general multivectors, rotors, quaternions, complex numbers, real numbers -- \emph{any subalgebra} of $\mathcal{G}(\mathbb{R}^n)$.

This is a very interesting result, accomplished due to the comprehensive analytic tools provided by Geometric Calculus. Recall that, in adaptive filtering theory, the transition from real-valued AFs to complex-valued AFs requires one to abide by the rules of differentiation with respect to complex-valued variables, represented by the \emph{Cauchy-Riemann conditions} (see~\cite[p.25]{Sayed08}). Similarly, quaternion-valued AFs require further differentiation rules that are captured by the \emph{Hamilton-real} (HR) calculus~\cite{Took_QLMS,Mandic2011,Jahanchahi2012} and its generalized version (GHR)~\cite{Dongpo2016}. Although those approaches are successful, each time the underlying algebra is changed, the analytic tools need an update as well. This is not the case if one resorts to GA and GC to address the minimization problem: the calculations are always performed the same way. 

\section{Mean-Square Analysis (Steady State)}\label{sec:mean-square-analysis}
The goal of the analysis is to derive an expression for the mean-square error (MSE) in steady-state of GAAFs via energy conservation relations (ECR)~\cite{Sayed08}. To achieve that, first some quantities and metrics are recast into GA.
\vspace*{-2mm} 

\subsection{Preliminaries}
\label{sec:data_model_GA}
A \emph{random multivector} is one whose blade coefficients are random variables. For instance, a random multivector in $\mathcal{G}(\mathbb{R}^3)$ is 
\eq
\small{
	\*A {=} \*a(0) {+} \*a(1)\gamma_1 {+} \*a(2)\gamma_2 {+} \*a(3)\gamma_3 {+} \*a(4)\gamma_{12} {+} \cdots {+} \*a(7) I,
}
\label{eq:random_multivector}
\qe   
where the terms $\*a(0), {\cdots}, \*a(7)$ are independent and identically distributed (i.i.d.) \emph{real-valued} random variables. By extension, \textit{random arrays} are arrays of random multivectors. 

The ECR technique is an energy balance in terms of the following (random) error quantities
\begin{equation}
\begin{dcases}
\Delta \*w_{i-1} \triangleq (w^o - \*w_{i-1})
\text{    weight-error array }\\
\boldsymbol{E}_a(i) = \rtb{u}_i \Delta w_{i-1}
\text{    \textit{a priori} estimation error }\\
\boldsymbol{E}_p(i) = \rtb{u}_i \Delta w_{i}
\text{    \textit{a posteriori} estimation error }
\end{dcases}
\label{eq:data_model_2}
\end{equation}
together with the AF's recursion.

The stationary data model is captured by the following set of assumptions
\begin{equation} 
\begin{array}{ll}
\text{\textbf{(1)} There exists an array of multivectors } w^{o} \text{ such that }\\ \hspace*{25mm}\boldsymbol{D}(i)=\rtb{u}_i w^{o} + \boldsymbol{V}(i)~\text{;}\\
\text{\textbf{(2)} The noise sequence } \lbrace\boldsymbol{V}(i)\rbrace \text{ is i.i.d. with constant}\\
\text{ variance } \expec\langle\b~{V}(i) \*{V}(i)\rangle;\\
\text{\textbf{(3)} The noise sequence } \lbrace\boldsymbol{V}(i)\rbrace \text{ is independent of } \boldsymbol{u}_j \text{ for all } \\i,j,~\text{and all other data;}\\
\text{\textbf{(4)} The initial condition } \*w_{-1} \text{ is independent of all }\\\lbrace\boldsymbol{D}(i),\boldsymbol{u}_i,\boldsymbol{V}(i)\rbrace~\text{;} \\
\text{\textbf{(5)} The regressor covariance matrix is } R_u = \expec\boldsymbol{u}_i\boldsymbol{u}^*_i > 0~\text{;}\\
\text{\textbf{(6)} The random quantities } \lbrace\boldsymbol{D}(i),\boldsymbol{u}_i,\boldsymbol{V}(i)\rbrace \text{ are zero mean.}
\end{array}
\label{eq:stationary_data_model} 
\end{equation}

As in linear algebra, the steady-state MSE in GA must be \emph{scalar-valued}. To this end, the MSE is defined as
\eq 
\text{MSE} = \xi \triangleq \lim_{i \to \infty} \expec\Big\langle\norm{\boldsymbol{E}(i)}^2\Big\rangle = \lim_{i \to \infty} \expec\Big\langle\b~{E}(i)\*{E}(i)\Big\rangle,
\label{eq:def_steadyState_MSE}
\qe 
where $\norm{\cdot}^2$, defined in~\eqref{eq:norm_of_array}, is applied to compactly denote the geometric product $\b~{E}(i)\*{E}(i)$.

From the stationary linear data model~\eqref{eq:stationary_data_model}, 
\vspace*{-2mm} 
\eq 
\eqarray 
\*E(i) \hspace*{-3mm}&{=} \*D(i) {-} \rtb{u}_i \iprod \*w_{i-1} {=} \rtb{u}_i\iprod(w^o {-} \*w_{i-1}) {+} \*V(i) {=} \*E_a(i) {+} \*V(i).
\label{eq:e_and_ea_and_v}
\qearray 
\vspace*{-1mm} 
\qe 
The term $\*E_a(i)$ is the \emph{a priori} error, from which the steady-state \emph{excess mean-square error} (EMSE) is defined, 
\vspace*{-1mm} 
\eq 
\text{EMSE} = \zeta \triangleq \lim_{i \to \infty} \expec\Big\langle\norm{\boldsymbol{E}_a(i)}^2\Big\rangle = \lim_{i \to \infty} \expec\Big\langle\b~{E}_a(i)\*{E}_a(i)\Big\rangle. 
\label{eq:def_steadyState_EMSE}
\vspace*{-1mm} 
\qe	 

From~\eqref{eq:e_and_ea_and_v}, $\expec\langle\b~{E}(i)\*{E}(i)\rangle = \expec\langle\b~{E}_a(i)\*{E}_a(i)\rangle + \expec\langle\b~{V}(i)\*{V}(i)\rangle$, since $\*V(i)$ is independent of any other random quantity and its samples are assumed to be drawn from a zero-mean white Gaussian process. Therefore, applying~\eqref{eq:def_steadyState_MSE} and \eqref{eq:def_steadyState_EMSE}, $\text{MSE} = \text{EMSE} + \expec\langle\b~{V}(i)\*{V}(i)\rangle$, analogous to the LA case. 

\subsection{Steady-State Analysis}
\label{sec:performance_GAAFs_1}

The ECR technique performs an interplay between the energies of the weight error array $\Delta\*w$ and the error $\*E$ at two successive time instants, $i-1$ (\emph{a priori}) and $i$ (\emph{a posteriori}). As a result, an expression for the \emph{variance relation} is obtained, which is then particularized for each AF of interest. For details on the ECR procedure, please refer to~\cite[p.228]{Sayed08}.    

Consider a GAAF whose update rule has the following general shape
\vspace*{-2mm}
\eq
\*w_{i} = \*w_{i-1} + \mu \*u_{i} f(\*E(i)),
\label{eq:ecr_1}
\vspace*{-2mm}
\qe
where $f(\cdot)$ is a multivector-valued function of the estimation error $\*E(i)$. Depending on the type of the GAAF (LMS, NLMS etc), $f(\cdot)$ assumes a specific value.         

Subtracting~\eqref{eq:ecr_1} from the optimal weight array $w^o$ yields
\vspace*{-2mm}
\eq
\Delta\*w_{i} = \Delta\*w_{i-1} - \mu \*u_{i} f(\*E(i)),
\label{eq:ecr_2}
\vspace*{-2mm}
\qe
in which $\Delta\*w_i = w^o - \*w_i$. Multiplying from the left by $\rtb{u}_{i}$ (array product) results in
\vspace*{-2mm}
\eq
\eqarray
\rtb u_{i} \iprod \Delta\*w_{i} &= \rtb u_{i} \iprod \left[\Delta\*w_{i-1} - \mu \*u_{i} f(\*E(i))\right]\\
\*E_p(i) &= \*E_a(i) - \mu \norm{\*u_{i}}^2 f(\*E(i)),
\label{eq:ecr_3}
\qearray
\vspace*{-2mm}
\qe
where $\*E_p(i) = \rtb u_{i} \iprod \Delta\*w_{i}$ is the \emph{a posteriori} error, $\*E_a(i) = \rtb u_{i} \iprod \Delta\*w_{i-1}$ is the \emph{a priori} error (See~\eqref{eq:data_model_2}), and in the last equation $\norm{\*u_{i}}^2 = \rtb u_{i} \iprod \*u_{i}$ (See~\eqref{eq:norm_of_array}).

The multivector $\norm{\*u_{i}}^2$ is assumed to be factorable into a product of invertible vectors\footnote{This assumes that $\norm{\*u_{i}}^2$ has only one type of grade: vector or bivector or trivector and so on. In practice, $\norm{\*u_{i}}^2$ can be a general multivector composed by different grades. However, such assumption allows for a clearer analysis procedure, and ultimately does not compromise the accuracy of the EMSE expression, as shown in the simulations.}~\cite[p. 14]{hestenes1987GAcalculus}, which guarantees the existence of a multiplicative inverse $\*\Gamma(i) \triangleq \big(\norm{\*u_{i}}^2\big)^{-1}$. This allows for solving~\eqref{eq:ecr_3} for $f(\*E(i))$
\vspace*{-2mm}
\eq
f(\*E(i)) = \mu^{-1} \*\Gamma(i) \left[\*E_a(i) - \*E_p(i)\right], 
\label{eq:ecr_4}
\vspace*{-2mm}
\qe
which substituted into~\eqref{eq:ecr_2} results in

\vspace*{-2mm}
\eq
\Delta\*w_{i} + \*u_{i}\*\Gamma(i)\*E_a(i) = \Delta\*w_{i-1} + \*u_{i}\*\Gamma(i)\*E_p(i).
\label{eq:ecr_6}
\vspace*{-2mm}
\qe

Taking the squared magnitude of both sides,
\vspace*{-1mm}
\eq
\underbrace{\left\lvert\Delta\*w_{i} + \*u_{i}\*\Gamma(i)\*E_a(i)\right\rvert^2}_{\text{LHS}~\eqref{eq:ecr_7}} = \underbrace{\left\lvert\Delta\*w_{i-1} + \*u_{i}\*\Gamma(i)\*E_p(i)\right\rvert^2}_{\text{RHS}~\eqref{eq:ecr_7}}.
\label{eq:ecr_7}
\vspace*{-1mm}
\qe

The left-hand side (LHS) is expanded as
\vspace*{-2mm}
\eq
\text{LHS} = \Big(\Delta\*w_{i} + \*u_{i}\*\Gamma(i)\*E_a(i)\Big)*\Big(\Delta\*w_{i} + \*u_{i}\*\Gamma(i)\*E_a(i)\widetilde{\Big)} 
\label{eq:ecr_8}
\vspace*{-2mm}
\qe
in which $*$ is the GA scalar product and \hspace*{1mm}$\t~{}$\hspace*{1mm} is the reverse. Further expansion gives
\vspace*{-2mm}
\eq
\eqarray
\hspace{-3mm}\text{LHS}&\hspace{-2mm}{=} \left\lvert\Delta\*w_{i}\right\rvert^2 {+}\left\lvert\*u_{i}\*\Gamma(i)\*E_a(i)\right\rvert^2\\ &\hspace{8mm}{+} \underbrace{\Delta\*w_{i}{*}\big(\b~E_a(i)\*\Gamma(i)\rtb u_{i}\big){+} \big(\*u_{i}\*\Gamma(i)\*E_a(i)\big){*} \t~\Delta \*w_{i}}_{\text{Sum of 3rd and 4th terms}}, 
\qearray
\label{eq:ecr_9}
\vspace*{-2mm}
\qe
in which $\*\Gamma(i) = \widetilde{\*\Gamma}(i)$ since $\t~{\norm{\*u_i}^2} = \norm{\*u_i}^2$ holds (See after~\eqref{eq:norm_of_array}). Applying the definition of GA scalar product and observing that the third and fourth terms of~\eqref{eq:ecr_9} are each other's reverse (i.e., their $0$-grades are the same), their sum can be written as $2\Big\langle\b~E_a(i)\*\Gamma(i)\rtb u_{i}\Delta\*w_{i}\Big\rangle$,

where the cyclic reordering property was used. Note that the term $\rtb u_{i}\Delta\*w_{i}$ is the definition of the \emph{a posteriori} error $\*E_p(i)$ (see~\eqref{eq:ecr_3}). This way, \eqref{eq:ecr_8} assumes the form
\vspace*{-2mm}
\eq
\left\lvert\Delta\*w_{i}\right\rvert^2 + 2\Big\langle\b~E_a(i)\*\Gamma(i)\*E_p(i)\Big\rangle {+} \left\lvert\*u_{i}\*\Gamma(i)\*E_a(i)\right\rvert^2. 
\label{eq:ecr_11}
\vspace*{-2mm}
\qe

A similar procedure allows for expanding the right-hand side (RHS) of \eqref{eq:ecr_7} as
\vspace*{-2mm}
\eq
\left\lvert\Delta\*w_{i-1}\right\rvert^2 + 2\Big\langle\b~E_p(i)\*\Gamma(i)\*E_a(i)\Big\rangle {+} \left\lvert\*u_{i}\*\Gamma(i)\*E_p(i)\right\rvert^2. 
\label{eq:ecr_12}
\vspace*{-1.5mm}
\qe

Substituting \eqref{eq:ecr_11} and \eqref{eq:ecr_12} into \eqref{eq:ecr_7}, and noting that the terms enclosed by the $0$-grade operator are each other's reverse (leading to mutual cancellation of their $0$-grades), 
\vspace*{-1.5mm}
\eq 
\left\lvert\Delta\*w_{i}\right\rvert^2 {+} \left\lvert\*u_{i}\*\Gamma(i)\*E_a(i)\right\rvert^2 {=} \left\lvert\Delta\*w_{i-1}\right\rvert^2 {+} \left\lvert\*u_{i}\*\Gamma(i)\*E_p(i)\right\rvert^2,  
\label{eq:ecr_13}
\vspace*{-1.5mm}
\qe
which is an energy relation balancing out \emph{a priori} and \emph{a posteriori} terms. Taking the expectation of the terms of~\eqref{eq:ecr_13} with respect to the random quantities $\*d(i)$ and $\*u_{i}$ results in
\eq
\expec\lvert\Delta\*w_{i}\rvert^2 + \expec\hspace{1mm}|\*u_{i}\*\Gamma(i)\*E_a(i)|^2 = \expec\lvert\Delta\*w_{i-1}\rvert^2 + \expec\hspace{1mm}|\*u_{i}\*\Gamma(i)\*E_p(i)|^2. 
\label{eq:ecr_16}
\qe   

Calculating the limit of~\eqref{eq:ecr_16} as $i \to \infty$ gives
\eq
\expec\hspace{1mm}|\*u_{i}\*\Gamma(i)\*E_a(i)|^2 = \expec\hspace{1mm}|\*u_{i}\*\Gamma(i)\*E_p(i)|^2, i\to\infty,
\label{eq:ecr_17}
\qe 
in which the steady-state condition $\expec\lvert\Delta\*w_{i}\rvert^2 = \expec\lvert\Delta\*w_{i-1}\rvert^2 = constant$ as $i\to\infty$ was employed. Plugging~\eqref{eq:ecr_3} into~\eqref{eq:ecr_17} results in
\eq
\expec\hspace{1mm}|\*u_{i}\*\Gamma(i)\*E_a(i)|^2 = \expec\hspace{1mm}|\*u_{i}\*\Gamma(i)(\*E_a(i) - \mu \norm{\*u_{i}}^2 f)|^2.
\label{eq:ecr_18}
\qe 

The right-hand side of~\eqref{eq:ecr_18} is expanded as
\eq
\expec\lvert\*u_{i}\*\Gamma(i)\*E_a(i)\rvert^2 - 2\mu\expec\Big\langle\*u_{i}\*\Gamma(i)\*E_a(i)\t~f\rtb u_{i}\Big\rangle + \mu^2\expec\lvert\*u_{i}f\rvert^2.
\normalsize
\label{eq:ecr_19}
\qe

Plugging~\eqref{eq:ecr_19} back into~\eqref{eq:ecr_18} and cancelling out the term $\expec\lvert\*u_{i}\*\Gamma(i)\*E_a(i)\rvert^2$ on both sides results in   
\eq
2\mu\expec\Big\langle\*u_{i}\*\Gamma(i)\*E_a(i)\t~f\rtb u_{i} \Big\rangle = \mu^2\expec\lvert\*u_{i}f\rvert^2.
\label{eq:ecr_21}
\qe  
Finally, applying the cyclic reordering on the left-hand side of~\eqref{eq:ecr_21} to make $\rtb u_{i}\*u_{i}\*\Gamma(i) = 1$ yields the \textit{variance relation}
\vspace*{-2mm}
\eq
\boxed{2\expec\Big\langle \*E_a(i)\t~f \Big\rangle = \mu\expec\big\lvert\*u_{i}f\big\rvert^2}.
\label{eq:ecr_22}
\vspace*{-1mm}
\qe  

For the GA-LMS, $f(\*E(i)) = \*E(i) = \*E_a(i) + \*V(i)$ (see~\eqref{eq:e_and_ea_and_v}). Substituting into~\eqref{eq:ecr_22},
\eq
2\expec\Big\langle \*E_a(i)\Big(\b~E_a(i) + \b~V(i)\Big)\Big\rangle = \mu\expec\Big\lvert\*u_{i}\Big(\*E_a(i) + \*V(i)\Big)\Big\rvert^2.
\label{eq:ecr_GA-LMS_1}
\qe 

The left-hand side of~\eqref{eq:ecr_GA-LMS_1} becomes
\vspace*{-2mm}
\eq
\eqarray
\hspace{-3mm}\text{LHS}~\eqref{eq:ecr_GA-LMS_1}&\hspace{-3mm}{=} 2\expec\big\langle\*E_a(i)\b~E_a(i)\big\rangle + 2\expec\big\langle \*E_a(i)\b~V(i)\big\rangle\\
&\hspace{-3mm}{=} 2\expec|\*E_a(i)|^2 + 2\left(\expec \*E_a(i) * \expec \b~V(i)\right) {=} 2\expec|\*E_a(i)|^2, 
\qearray
\label{eq:ecr_GA-LMS_2}
\vspace*{-1mm}
\qe 
where once more the independence of $\*V(i)$ was utilized, and $\expec \b~V(i) = \expec \*V(i)= 0$ (its entries are drawn from a zero-mean white Gaussian process).

The right-hand side of~\eqref{eq:ecr_GA-LMS_1} is expanded as
\vspace*{-2mm}
\eq
\eqarray
\hspace{-2mm}\text{RHS}~\eqref{eq:ecr_GA-LMS_1}&\hspace{-3mm}= \mu\expec\left[\*u_i(\*E_a(i) + \*V(i)) * (\b~E_a(i) + \b~V(i))\rtb u_i\right]{=}\\
&\hspace*{-3mm}=\mu\expec\big\langle\norm{\*u_{i}}^2\norm{\*E_a(i)}^2\big\rangle + \mu\expec\big\langle \norm{\*u_{i}}^2\norm{\*V(i)}^2\big\rangle.
\qearray
\label{eq:ecr_GA-LMS_3}
\vspace*{-1mm}
\qe 

Substituting~\eqref{eq:ecr_GA-LMS_2} and~\eqref{eq:ecr_GA-LMS_3} into~\eqref{eq:ecr_GA-LMS_1} yields
\eq
2\expec|\*E_a(i)|^2 = \mu\expec\big\langle\norm{\*u_{i}}^2\norm{\*E_a(i)}^2\big\rangle + \mu\expec\big\langle \norm{\*u_{i}}^2\norm{\*V(i)}^2\big\rangle.
\label{eq:ecr_GA-LMS_5}
\qe 
Observing that $\expec|\*E_a(i)|^2 {=} \expec\big\langle\norm{\*E_a(i)}^2\big\rangle$, \eqref{eq:ecr_GA-LMS_5} is rewritten as
\eq
\expec\big\langle(2 - \mu\norm{\*u_{i}}^2)\norm{\*E_a(i)}^2\big\rangle = \mu\expec\big\langle \norm{\*u_{i}}^2\norm{\*V(i)}^2\big\rangle.
\label{eq:ecr_GA-LMS_6}
\qe 
Adopting the \emph{separation principle} (see~\cite[p.245]{Sayed08}), i.e., in steady state $\norm{\*u_{i}}^2$ is independent of $\*E(i)$ (and consequently of $\*E_a(i)$), \eqref{eq:ecr_GA-LMS_6} becomes
\eq
\big\langle(2 - \mu\expec\norm{\*u_{i}}^2)\expec\norm{\*E_a(i)}^2\big\rangle = \mu\big\langle \expec\norm{\*u_{i}}^2\expec\norm{\*V(i)}^2\big\rangle.
\label{eq:ecr_GA-LMS_7}
\qe 

From the Appendix it is known that the quantities $\expec\norm{\*u_{i}}^2$ and $\expec\norm{\*V(i)}^2$ depend of the underlying (sub)algebra. For $\mathcal{G}(\mathbb{R}^n)$, $\expec\norm{\*V(i)}^2$ and $\expec\norm{\*u_{i}}^2$ are obtained via~\eqref{eq:appendix_expec_vv} and~\eqref{eq:expec_uu_4} respectively, which substituted into~\eqref{eq:ecr_GA-LMS_7} yields  
\eq
(2 - \mu M(2^n\sigma^2_u))\big\langle\expec\norm{\*E_a(i)}^2\big\rangle = \mu M(2^n\sigma^2_u)(2^n\sigma^2_v),
\label{eq:ecr_GA-LMS_8}
\qe 
where $M$ is the regressor length (number of filter taps), and $\sigma^2_u, \sigma^2_v$ are  It is important to notice that since~\eqref{eq:expec_uu_4} is obtained considering the coefficients of the regressor entries are drawn from a circular Gaussian process (see appendix), the present analysis holds only for that kind of input.   

Finally, the expression for the GA-LMS steady-state EMSE using the complete algebra $\mathcal{G}(\mathbb{R}^n)$ is given by
\vspace*{-2mm}
\eq
\boxed{\underset{\hbox{\scriptsize LMS}}{\mathlarger{\zeta}} \hspace{2mm} = \hspace{2mm} \dfrac{\mu M 4^n\sigma^2_u\sigma^2_v}{2 - \mu M 2^n\sigma^2_u}, \hspace{2mm}i\to \infty}.
\label{eq:ecr_GA-LMS_9}
\vspace*{-1mm}
\qe

Table~\ref{tab:EMSE_list_GA-LMS} summarizes the EMSE for the even subalgebras of interest. Notice that for $\mathcal{G}^{+}(\mathbb{R})$ the EMSE for the LMS with real-valued entries is recovered (compare to Eq. 16.10 in~\cite[p.246]{Sayed08} for white Gaussian inputs). To obtain the respective MSE, one should add $\expec\langle\b~V(i)\*V(i)\rangle {=} \expec\norm{\*V(i)}^2$ to the EMSE value, as aforementioned. 

\begin{table}[]
	\centering
	\caption{Steady-state EMSE for the GA-LMS using even algebras.}
	\label{tab:EMSE_list_GA-LMS}
	\begin{tabular}{ccc}
		\hline
		\multicolumn{1}{l}{$\mathcal{G}^{+}(\mathbb{R}^n)$ (Even Algebras)}                  & \multicolumn{2}{l}{$\dfrac{\mu M \Big[\sum_{k}^{}\binom n{2k}\Big]^2 \sigma^2_u\sigma^2_v}{2 - \mu M \sigma^2_u\sum_{k}^{}\binom n{2k}},\text{ for } k \in \mathbb{N}$ \vspace*{1mm}} \\ \hline
		\multicolumn{1}{c|}{$\mathcal{G}^{+}(\mathbb{R}^3)$ (Quaternions) } & \multicolumn{1}{c|}{$\mathcal{G}^{+}(\mathbb{R}^2)$ (Complex)}                                          & $\mathcal{G}^{+}(\mathbb{R})$ (Real)                         \\
		\multicolumn{1}{c|}{$\dfrac{16 \mu M \sigma^2_u\sigma^2_v}{2 - 4 \mu M \sigma^2_u}$} & \multicolumn{1}{c|}{$\dfrac{4 \mu M \sigma^2_u\sigma^2_v}{2 - 2\mu M \sigma^2_u}$}                      & $\dfrac{\mu M\sigma^2_u\sigma^2_v}{2 - \mu M\sigma^2_u}$                     
	\end{tabular}
\end{table}

\section{Simulations}
\label{cap:applications_GAAFs}
This Section shows the performance of the computational implementation of GAAFs in a system identification task~\footnote{All the AFs were implemented using the \emph{Geometric Algebra ALgorithms Expression Templates} (GAALET)~\cite{seybold_gaalet_2010}, a C++ library for evaluation of GA expressions, and OpenGA~\cite{OpenGA}. The reader is encouraged to follow the instructions on \href{https://openga.org/ieeetsp.html}{openga.org/ieeetsp.html} in order to reproduce the simulations.}. The \emph{optimal weight array} $w^o$ to be estimated has $M$ multivector-valued entries (number of taps), namely $W^o_j$, $j = 1,\cdots,M$,	$w^o = \begin{bmatrix} W^o_1 \hspace*{1mm} W^o_2 \cdots W^o_M \end{bmatrix}^{T}$. Each case studied in the sequel (multivector, rotor, complex, and real entries) adopts a different value for $W^o_j$ (highlighted ahead).

As aforementioned, the measurement noise multivector $\*V$ has each of its coefficients drawn from a white Gaussian stationary process with variance $\sigma^2_V$.  

\subsection{Multivector Entries}
\label{ssec:app_GAAFs_1_MV_entries}
The underlying geometric algebra in this case is $\mathcal{G}(\mathbb{R}^{n})$, with $n=3$, i.e., the one whose multivectors are described by basis~\eqref{eq:orthonormal_basis_R3}. The derivation of the GAAFs puts no restriction on the values the vector space dimension $n$ can assume. However, the case $n=3$ (generating a GA with dimension $8$) is an example that captures the core idea of this work: the GAAFs can estimate hypercomplex quantities which generalize real, complex, and quaternion entries.    

In this case all the multivector entries of $w^o$ are the same, namely $W^o_j = \{0.55 + 0\gamma_1 + 1\gamma_2 + 2\gamma_3 + 0.71\gamma_{12} + 1.3\gamma_{23} + 4.5\gamma_{31} + 3I\}$ for $j=1,\cdots,M$. Those values were selected in an aleatory manner. 
Note that the coefficient of $\gamma_1$ is zero. However, it is displayed to emphasize the structure of the $\mathcal{G}(\mathbb{R}^{3})$ basis.    

\begin{figure}[t!]
	\centering
	\includegraphics[width=0.4\textwidth]{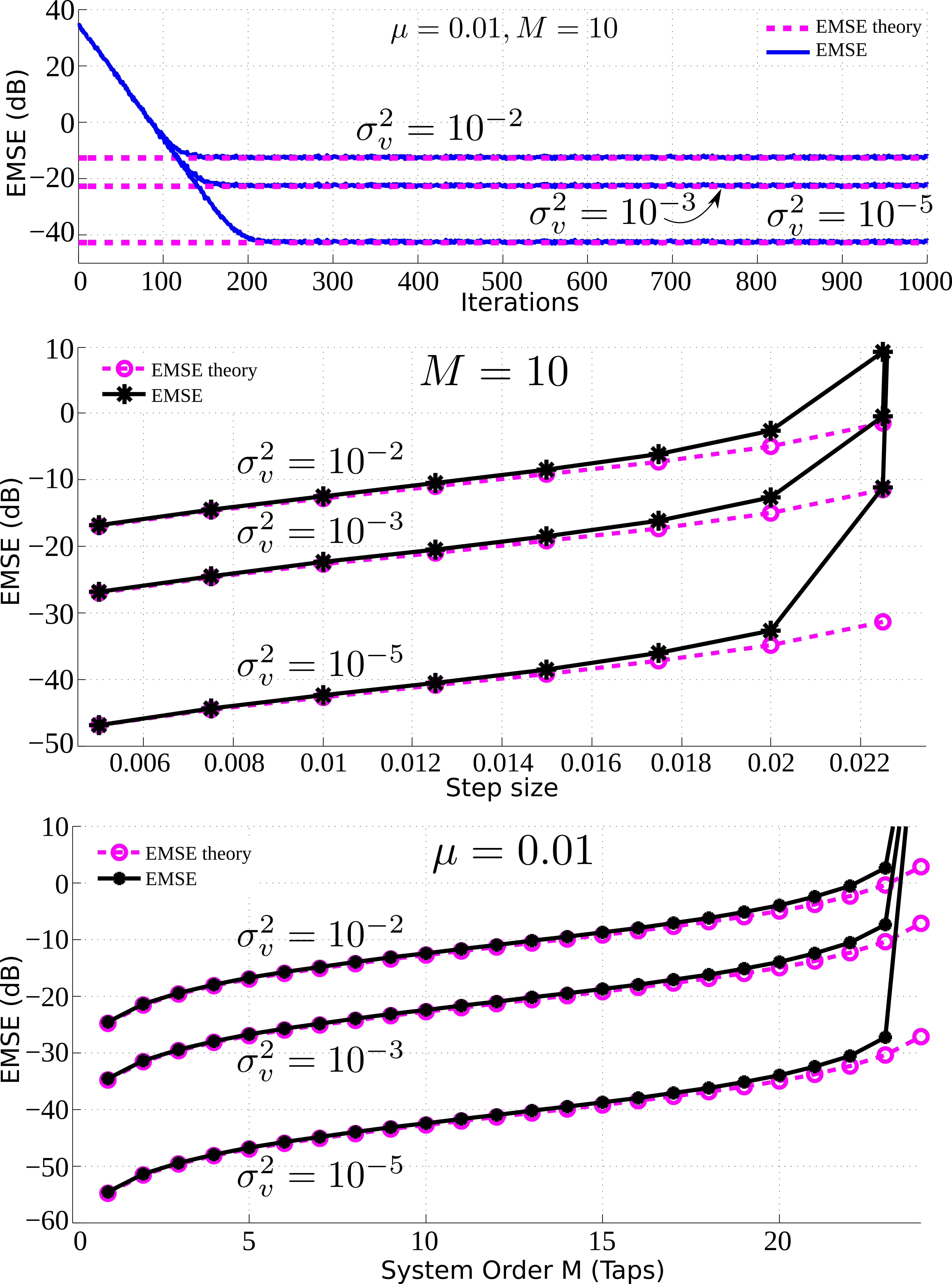}
	\caption{GA-LMS for multivector entries. (Top) EMSE learning curves for $\mu = 0.01$, $M=10$, and $\sigma^2_V = \{10^{-2},10^{-3},10^{-5}\}$ averaged over $100$ experiments. (Middle) Steady-state EMSE as a function of the step size $\mu$ for $M=10$. Note the AF stands on the brink of divergence for $\mu > 0.02$. (Bottom) Steady-state EMSE as a function of the taps $M$ for $\mu=0.01$. For $M > 23$ the AF starts diverging. The simulated steady-state value is obtained by averaging the last $200$ points of the ensemble-average learning curve for each $\mu$ and $M$.}
	\label{fig:GA-LMS_learningCurves_mv}
\end{figure}

\figurename~\ref{fig:GA-LMS_learningCurves_mv}~(top) shows several learning curves (MSE and EMSE) for the GA-LMS estimating the weight array $w^o$ with $M = 10$. The step size value is $\mu = 0.005$ for all simulated curves. Notice the perfect agreement between the theoretical error levels (obtained from~\eqref{eq:ecr_GA-LMS_9} with $n=3$) and the simulated steady-state error. Those experiments show that the GA-LMS is indeed \emph{capable of estimating multivector-valued quantities}, supporting what was previously devised in Section~\ref{cap:GAAFs_standard}. \figurename~\ref{fig:GA-LMS_learningCurves_mv}~(bottom) depicts the steady-state error as a function of the system order (number of taps) $M$ for $\sigma^2_V = \{10^{-2},10^{-3},10^{-5}\}$. Theory and experiment agree throughout the entire tested range $M = [1, 40]$.   

\subsection{Rotor, Complex, and Real Entries}
\label{ssec:app_GAAFs_1_Rotor_entries}
\begin{figure}[t!]
	\centering
	\includegraphics[width=0.18\textwidth, angle=-90]{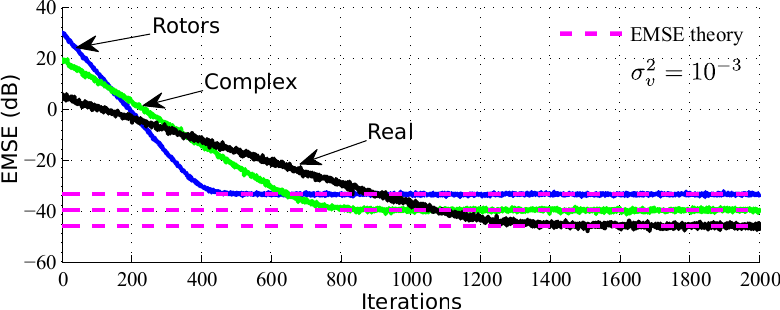}
	\caption[GA-LMS (rotor entries)]{EMSE learning curves for $M=10$, $\mu = 0.005$, and $\sigma^2_V = 10^{-3}$ for rotor, complex, and real entries ($100$ experiments). Notice how the theoretical and experimental values agree.}
	\label{fig:MSE_and_EMSE_GALMS_rotor}
\end{figure}

\figurename~\ref{fig:MSE_and_EMSE_GALMS_rotor} depicts the EMSE curves for three types of GA-LMS: $\mathcal{G}^{+}(\mathbb{R}^{3})$ (isomorphic to quaternions), $\mathcal{G}^{+}(\mathbb{R}^{2})$ (isomorphic to complex numbers), and $\mathcal{G}^{+}(\mathbb{R})$ (isomorphic to reals). All filters have $M = 10$, $\mu = 0.005$, and $\sigma^2_V = 10^{-3}$. However, each AF has a different entries for the optimal weight array $w^o$: $W^o_j = \{0.55 + 0.71\gamma_{12} + 1.3\gamma_{23} + 4.5\gamma_{31}\}$ for rotors; $W^o_j = \{0.55 + 0.71\gamma_{12}\}$ for complex; and $W^o_j = \{0.55\}$ for reals, with $j=1,\cdots,M$. 

The AFs are shown to be capable of estimating their respective optimal weight arrays with very good agreement with the theoretical value (Table~\ref{tab:EMSE_list_GA-LMS}). Due to the aforementioned isomorphism, those filters become alternatives to their LA counterparts, namely quaternion-LMS (QLMS)~\cite{Took_QLMS,Mandic2011,Jahanchahi2012,Dongpo2016}, complex-LMS (CLMS)~\cite{Widrow1975}, and real-LMS (LMS)~\cite{Sayed08,Diniz2013}.  

\section{Conclusion} 
\label{cap:conclusion}

The formulation of GA-based adaptive techniques is still in its infancy. The majority of AF algorithms available in the literature resorts to specific subalgebras of GA (real, complex numbers and quaternions). Each of them requires an specific set of tools in order to pose the problem and perform calculus. In this sense, the development of the GAAFs is an attempt to unify those different adaptive-filtering approaches under the same mathematical language. Additionally, as shown throughout the text, GAAFs have improved estimation capabilities since they are not limited to $1$-vector estimation (like LA-based AFs). Instead, they can naturally estimate any kind of multivector. Also, for the GA-LMS, the shape of its update rule is invariant with respect to the multivector subalgebra. This is only possible due to the use of GA and GC.

On top of the theoretical contributions, the experimental validation provided in Section~\ref{cap:applications_GAAFs} shows that the GAAFs are successful in a system identification task. Nevertheless, it is expected that any estimation problem posed in terms of hypercomplex quantities will benefit from this work. For instance, GAAFs may be useful in data fusion, where different signals are supposed to be integrated in the same ``package'' and then processed. The multivector (and by extension the array of multivectors) can be interpreted as a fundamental information package that aggregates scalar, vector, bivector, and so on, quantities.  

New types of GAAFs are currently under study, particularly the NLMS and RLS variants. Together with the transient analysis of GA-LMS and the introduction of noncircularity conditions~\cite{MandicBook,Neto_WLQLMS_2011}, they figure as subjects of future publications. Also, given the connection between exterior and tensor algebras~\cite{yokonuma1992tensor}, it would be interesting to investigate how GAAFs and tensor-product AFs~\cite{Rupp2015_Icassp,Rupp2015_Eusipco} are related.

\appendix%
\emph{Calculating the expectation} $\expec\norm{\*V(i)}^2$: Take a random multivector $\*V \in \mathcal{G}(\mathbb{R}^3)$ (see~\eqref{eq:random_multivector}) $\*V = \*v(0) + \*v(1)\gamma_1 + \cdots + \*v(6)\gamma_{31} + \*v(7) I$, where $\*v(k)$, $k=0,\cdots,7$, are i.i.d. real-valued random variables drawn from a zero-mean and stationary white Gaussian process. Performing the geometric product $\b~{V}\*{V} = \norm{\*V}^2$ and calculating its expectation results in	
\vspace*{-2mm} 
{\small
	\eq 
	\expec \b~{V}\*{V} = \expec\*v^2(0) + \expec\*v^2(1) + \expec\*v^2(2) + \expec\*v^2(3) + \cdots + \expec\*v^2(7), 
	\label{eq:expec_vv}
	\vspace*{-1mm} 
	\qe 
}
in which the expectations of the cross-terms are zero due to the i.i.d. assumption above. Each term $\expec\*v^2(k)$, is said to be the variance of $\*v(k)$ and denoted $\expec\*v^2(k) \triangleq \sigma^2_v$. This way, \eqref{eq:expec_vv} becomes $\expec \b~{V}\*{V} = 8\sigma^2_v, \hspace{1mm} \*V \in \mathcal{G}(\mathbb{R}^3)$. Note that in general $\expec\b~{V}\*{V} = \text{dim}\{\mathcal{G}^g(\mathbb{R}^n)\}\sigma^2_v$ for $\*V \in \mathcal{G}^g(\mathbb{R}^n)$, in which $\mathcal{G}^g(\mathbb{R}^n)$ can be \emph{any subspace} of $\mathcal{G}(\mathbb{R}^n)$. When the complete geometric algebra $\mathcal{G}(\mathbb{R}^n)$ is used,
\vspace*{-2mm}
{\small
	\eq
	\boxed{\expec\norm{\*V(i)}^2 = \expec\b~{V}\*{V} = 2^n\sigma^2_v, \hspace{1mm} \*V \in \mathcal{G}(\mathbb{R}^n)}.
	\label{eq:appendix_expec_vv}
	\vspace*{-1mm} 
	\qe 
}

\emph{Calculating the expectation} $\expec\norm{\*u_i}^2$: The regressor array $\*u_i$ is a collection of $M$ random multivectors $\*U_j \in \mathcal{G}(\mathbb{R}^n), j=1,\cdots,M$. Analogously to the LA case, the regressor covariance matrix is calculated as $R_u = \expec\boldsymbol{u}_i\boldsymbol{u}^*_i$. Its trace is $Tr(R_u) = \expec\rtb{u}_i\*u_i = \expec\norm{\*u_i}^2$, a multivector-valued quantity obtained via~\eqref{eq:array_product}, $\expec \rtb{u}\*{u} = \expec\b~{U}_1\*{U}_1 + \expec\b~{U}_2\*{U}_2 + \cdots + \expec\b~{U}_M\*{U}_M.$

For the special case $\*U_j \in \mathcal{G}(\mathbb{R}^3)$, the geometric product $\b~{U}_j\*{U}_j$ is  
{\scriptsize
	\eq 
	\eqarray 
	\b~{U}_j\*{U}_j =\hspace*{-3mm} 
	&\*u^2(j,0) + \*u(j,0)\*u(j,1)\gamma_1 + \cdots + \*u(j,0)\*u(j,7)I +\\
	&\*u^2(j,1) + \*u(j,0)\*u(j,1)\gamma_1 + \cdots + \*u(j,1)\*u(j,5)I +\\
	&\hspace{7mm}\vdots\hspace{17.5mm}\vdots\hspace{22.5mm}\vdots\\
	&\*u^2(j,7) + \*u(j,7)\*u(j,5)\gamma_1 + \cdots + \*u(j,7)\*u(j,0)I, \\
	\qearray
	\label{eq:UU}
	\qe  		
}%
where each real coefficient $\*u(j,k), k=0,\cdots,7$, is an i.i.d. random variable drawn from a zero-mean and stationary white Gaussian process. Thus, $\expec\b~{U}_j\*{U}_j \hspace*{-1mm}= \expec\*u^2(j,0) + \expec\*u^2(j,1) + \expec\*u^2(j,2) + \cdots + \expec\*u^2(j,7),$

since the expectations of the cross-terms in~\eqref{eq:UU} are zero. Each term $\expec\*u^2(j,k)$ is said to be the variance of $\*u(j,k)$ and denoted $\expec\*u^2(j,k) \triangleq \sigma^2_u$ (regressors are assumed to have shift structure). Note that this result is also obtained if $\*{U}_j$, $j = 1,\cdots,M$, is a circular Gaussian random multivector (the circularity condition~\cite[p. 8]{Sayed08} for complex-valued random variables is extended here to encompass random multivectors.). Such case considers that the coefficients of a random multivector are independent Gaussian random variables. This way, $\expec \b~{U}_j\*{U}_j = 8\sigma^2_u$, which yields $\expec \rtb{u}\*u = M (8\sigma^2_u), \hspace*{1mm} \*U_j \in \mathcal{G}(\mathbb{R}^3)$. Note that in general $\expec\rtb{u}\*{u} = M(\text{dim}\{\mathcal{G}^g(\mathbb{R}^n)\}\sigma^2_u)$ for $\*U_j$ belonging to $\mathcal{G}^g(\mathbb{R}^n)$ (any subspace of $\mathcal{G}(\mathbb{R}^n)$). When $\*U_j \in \mathcal{G}(\mathbb{R}^n)$, 
\vspace*{-2mm} 
\eq 
\boxed{\expec\norm{\*u}^2 = \expec\rtb{u}\*{u} = M(2^n\sigma^2_u), \hspace*{1mm} \*U_j \in \mathcal{G}(\mathbb{R}^n)}. 
\label{eq:expec_uu_4}
\vspace*{-2mm} 
\qe

\ifCLASSOPTIONcaptionsoff
  \newpage
\fi



%
\bibliographystyle{IEEEbib}
\bibliography{refs}
\vspace*{-1cm}

\end{document}